

 \magnification=\magstephalf  \hoffset=0.5truein 

 \input amstex \documentstyle{amsppt}            

\refstyle{A}

\redefine\Re{\operatorname{Re}}\redefine\Im{\operatorname{Im}}

\define\cI{{\Cal I}}\define\cJ{{\Cal J}}\define\cK{{\Cal K}}

\define\bbR{{\Bbb R}}\define\bbC{{\Bbb C}}
\define\bbZ{{\Bbb Z}}\define\bbP{{\Bbb P}}

\define\ebold{{\bold e}}

\define\euso{\operatorname{\frak{so}}}\define\eusu{\operatorname{\frak{su}}}

\define\eugl{\operatorname{\frak{gl}}}\define\eug{\operatorname{\frak g}}
\define\euh{\operatorname{\frak h}}\define\eus{\operatorname{\frak s}}
\define\SO{\operatorname{SO}}\define\SL{\operatorname{SL}}
\define\SU{\operatorname{SU}}\define\Un{\operatorname{U}}
\define\GL{\operatorname{GL}}\define\Symp{\operatorname{Sp}}
\define\tr{\operatorname{tr}}\define\Or{\operatorname{O}}
\define\codim{\operatorname{codim}}
\define\Vol{\operatorname{Vol}}\define\Gr{\operatorname{Gr}}
\define\G{\operatorname{G}}\define\Spin{\operatorname{Spin}}

\define\lhk{\mathbin{\hbox{\vrule height1.4pt width4pt depth-1pt 
             \vrule height4pt width0.4pt depth-1pt}}}
\define\w{{\mathchoice{\,{\scriptstyle\wedge}\,}{{\scriptstyle\wedge}}
      {{\scriptscriptstyle\wedge}}{{\scriptscriptstyle\wedge}}}}

\topmatter
\title Calibrated Embeddings\\
in the Special Lagrangian\\
and Coassociative Cases
\endtitle
\author
Robert L. Bryant
\endauthor
\affil
Duke University
\endaffil
\address
Department of Mathematics,
Duke University,
Box 90320,
Durham, NC 27708-0320
\endaddress
\date
December 27, 1999
\enddate

\email
bryant\@math.duke.edu
\endemail



\dedicatory
Dedicated to the memory of Alfred Gray
\enddedicatory

\thanks 
The research for this article was made possible by support 
from the National Science Foundation through grant DMS-9870164
and from Duke University.
\endthanks

\abstract
Every closed, oriented, real analytic Riemannian 
$3$-manifold can be isometrically embedded as a special Lagrangian 
submanifold of a Calabi-Yau $3$@-fold, even as the real locus of an 
antiholomorphic, isometric involution.  Every 
closed, oriented, real analytic Riemannian $4$-manifold whose bundle of 
self-dual $2$@-forms is trivial can be isometrically embedded as a 
coassociative submanifold in a $\G_2$-manifold, even as the fixed locus
of an anti-$\G_2$ involution.

These results, when coupled with McLean's analysis of the moduli spaces 
of such calibrated submanifolds, yield a plentiful supply of examples of 
compact calibrated submanifolds with nontrivial deformation spaces.
\endabstract

\rightheadtext{Calibrated Embedding}

\endtopmatter

\document

\head 0. Introduction  \endhead

\subhead 0.0. Calibrations \endsubhead 
Let~$(M,g)$ be a Riemannian manifold and let~$\phi$ be a
{\it calibration\/} on~$M$, i.e., a closed $p$@-form~$\phi$ 
with the property that
$$
\phi(\ebold_1,\dots,\ebold_p)\le 1 
\tag1
$$
for all $p$-tuples of orthonormal vectors~$\ebold_1,\dots,\ebold_p\in T_xM$
for any~$x\in M$.  

Any oriented~$p$@-plane~$E\subset T_xM$, 
satisfies~$\iota^*_E\phi\le\Omega_E$, 
where~$\iota_E:E\to T_xM$ is the inclusion and~$\Omega_E$ 
is the canonical oriented unit volume element in~$\Lambda^n(E^*)$
induced by the metric~$g$ and the given orientation of~$E$. 
An oriented $p$@-plane~$E\subset T_xM$ is said to be
{\it calibrated\/} by~$\phi$ if~$\iota_E^*\phi=\Omega_E$.  

For any oriented $p$@-dimensional submanifold~$N\subset M$, the inequality~ 
$\iota_N^*\phi\le\Omega_N$ holds, where~$\Omega_N$ is the induced oriented 
unit volume form on~$N$ and~$\iota_N:N\to M$ is the inclusion mapping.  
The submanifold~$N$ is said to be {\it calibrated\/} by~$\phi$ 
if each of its oriented tangent planes is, i.e., if $\iota_N^*\phi=\Omega_N$.

If $N$ is a closed, oriented $p$@-dimensional submanifold of~$M$ 
that is calibrated by~$\phi$, then~$N$ is 
{\it homologically absolutely minimizing}:  
If~$N'$ is any closed oriented 
$p$-dimensional submanifold homologous to~$N$
then~$\Vol(N)\le\Vol(N')$, with equality if and only if~$N'$ 
is also calibrated by~$\phi$.  

This follows by writing~$N'-N = \partial C$ for some 
(oriented) $(p{+}1)$@-cycle~$C$ and then observing that
$$
\align
\Vol(N') &= \int_{N'} \Omega_{N'} \ge \int_{N'} \phi
                \qquad \text{(since $\phi$ is a calibration),}\\
         &= \int_{N} \phi + \int_C d\phi
                \qquad \text{(by Stokes' Theorem),}\\
         &= \int_{N} \Omega_N + 0
     \qquad \text{(since $N$ is calibrated by~$\phi$ and since $d\phi=0$),}\\
		 &= \Vol(N)\,.\\
\endalign
$$
Since equality would imply that~$\Omega_{N'}=\iota_{N'}^*\phi$ almost 
everywhere, it would also imply that~$N'$ must be calibrated by~$\phi$.  

Thus, calibrations provide a method of proving that certain submanifolds of
Riemannian manifolds are not only minimal, but are actually homologically
minimizing, a much stronger condition.  The theory of calibrations was
developed extensively by Harvey and Lawson~\cite{HL} and is a fundamental
tool in the theory of minimal submanifolds.

Although the `generic' calibration~$\phi$ will not calibrate any 
submanifolds, Harvey and Lawson show that there are many calibrations
for which calibrated submanifolds are plentiful, at least
locally.  However, even in these cases, the study of compact or 
complete calibrated submanifolds is usually less complete.

The best-known example is that of a K\"ahler form and
its divided powers on a complex manifold (cf.~\S0.1), where 
the calibrated submanifolds are simply the complex submanifolds.
Of course, there are many complex submanifolds of a given complex
manifold, at least locally.  Compact complex submanifolds
are much more rigid.  However, since, in this case, the moduli space of 
compact calibrated manifolds can be studied using techniques from complex 
geometry, our knowledge about these moduli spaces is rather extensive.  
In particular, vanishing theorems in the holomorphic category can be 
brought to bear to show that the moduli space of compact complex 
submanifolds is smooth in many cases.

The work of McLean~\cite{Mc} has provided two more
examples in which the local geometry of the moduli space of
compact calibrated submanifolds can be shown to be smooth.  They
are the special Lagrangian calibration~(cf.~\S0.2) and the 
coassociative calibration~(cf.~\S0.3).  These cases have 
been extensively explored in the literature, partly for
their intrinsic interest~(\cite{Br1}, \cite{Br3}, \cite{Hi1},
\cite{Hi2}) and partly for their interest for string theory
and mirror symmetry (\cite{Ac1, Ac2}, \cite{BSh}, \cite{Li}, 
\cite{SYZ}, to give just a sample). 

However, explicit examples of compact calibrated submanifolds in these
cases are not so easy to construct, which makes the theory difficult
to investigate.  It is not even clear what sorts of restrictions there
might be on the intrinsic geometry of calibrated submanifolds
in these cases.  

Part of the difficulty is that it is not so easy to construct nontrivial 
examples of such calibrations in the first place. 
For example, the special Lagrangian calibration is defined on
a Calabi-Yau manifold and non-trivial compact examples of these latter
spaces were not known to exist until the celebrated work 
of Yau~\cite{Ya}.%
\footnote{applied to string theory in~\cite{CHSW}}
The coassociative calibration is defined on
a Riemannian 7-manifold with holonomy~$\G_2$ and nontrivial local examples
of such metrics were not known to exist until much later~\cite{Br2}.
Compact examples came later still, with the work of Joyce~\cite{Jo}.%
\footnote{applied to string theory in~\cite{PT} and~\cite{SV}}  
In either case, compact examples of these calibrated
submanifolds are still very rare.

In this article, I will show that there are many nontrivial
examples of compact calibrated submanifolds~$N$ in these two cases.
However, the ambient manifold~$M$ that will be constructed will not
generally be compact or complete.  Instead, it will be more like
a `germ' of a Calabi-Yau or $\G_2$@-manifold that forms a neighborhood
of the given~$N$.  Specifically, I show that any closed oriented, 
real analytic Riemannian $3$-manifold can be isometrically embedded as 
a special Lagrangian submanifold of some Calabi-Yau $3$-fold (Theorem~1) 
and that any closed oriented, real analytic Riemannian 4-manifold whose 
bundle of self-dual 2-forms is trivial can be isometrically
embedded as a coassociative submanifold of a $\G_2$@-manifold (Theorem 2).

While the errors and omissions in this article are my own, I would like 
to thank Dave Morrison for his valuable advice on references to 
the literature.

\subhead 0.1. The almost K\"ahler case \endsubhead 
The most familiar example of a nontrivial calibration is to be found
in K\"ahler geometry.  Let~$(M,g,\Omega)$ be an almost-K\"ahler manifold. 
I.e., $g$ is a Riemannian metric on~$M^{2m}$ and~$\Omega$ is a closed 2-form 
on~$M$ with the property that the skewsymmetric endomorphism~$J:TM\to TM$ 
that satisfies $g(v,w) = \Omega(v,Jw)$ for all~$v,w\in T_xM$ and all~$x\in M$
defines a complex structure on~$TM$.  

The well-known Wirtinger inequality~\cite{HL} implies that
for~$p\le m$ the $2p$-form
$$
\Omega^{[p]} = \frac1{p!}\,\Omega^p
$$
is a calibration on~$M$ and that an oriented $2p$-plane~$E\subset T_xM$ is 
calibrated by~$\Omega^{[p]}$ if and only if it is a complex $p$-dimensional
subspace of~$T_xM$ endowed with its natural orientation as a 
complex vector space.  Thus, the $\Omega^{[p]}$@-calibrated submanifolds
are the almost-complex%
\footnote{Also known as `pseudo-holomorphic'.}
$p$-manifolds in~$M$.  

In the most interesting case, when~$J$ is an integrable complex structure
(i.e., the K\"ahler case),
the $p$-dimensional complex submanifolds of~$M$ are plentiful, at least
locally.  Moreover, given a compact $p$-dimensional complex submanifold
$N^{2p}\subset M$, the moduli space of `nearby' $p$-dimensional complex 
submanifolds is an analytic variety and can be described purely 
in terms of the underlying complex geometry.

Note also that, in the general almost-complex case,
an almost-complex $N^{2p}\subset M$ inherits its own 
almost-K\"ahler structure by pullback from~$M$.  This structure will
be K\"ahler whenever the ambient structure is K\"ahler.

\subhead 0.2. The almost special Lagrangian case \endsubhead 
Let
$$
\align
g_0 &= dz^1\circ d{\bar z}^1 + \cdots +dz^m \circ d{\bar z}^m \\
\omega_0 &= \frac\imath2
      \bigl(dz^1\w d{\bar z}^1 + \cdots +dz^m\w d{\bar z}^m \bigr)\\
\Upsilon_0 &= dz^1\w \cdots dz^m\\
\endalign
$$
be the standard metric, K\"ahler form, and holomorphic volume form 
on~$\bbC^m$.  The subgroup of~$\GL(m,\bbC)$ that preserves these forms 
is~$\SU(m)$.  For later use, I want to make an explicit identification
of~$\bbC^m$ with~$\bbR^{2m}$, namely~$z = x+\imath\,y\in\bbC^m$ with
$x,y\in\bbR^m$ will be identified with~${x\choose y}\in\bbR^{2m}$.

\subsubhead 0.2.1. The calibration \endsubsubhead
Harvey and Lawson~\cite{HL} show that the $m$@-form~$\phi_0=\Re\Upsilon_0$ 
is a calibration on~$\bbC^{2m}$ and that, moreover, an 
$m$@-plane~$E\subset T_0\bbC^m\simeq\bbC^m$ is $\phi_0$@-calibrated if 
and only if there exists an~$A\in\SU(m)$ so that~$A(E)=\bbR^m\subset\bbC^m$.
Thus, any~$\phi_0$@-calibrated~$E$ satisfies~$\iota_E^*(\omega_0)=0$, 
i.e.,~$E$ is an $\omega_0$@-Lagrangian $m$@-plane.  

Harvey and Lawson further show that any $\omega_0$@-Lagrangian~$m$@-plane~$E$ 
satisfies~$\iota_E^*(\Upsilon_0)=\lambda(E)\,\Omega_E$ for some complex 
number~$\lambda(E)$ satisfying~$|\lambda(E)|=1$.  Thus, 
an~$\omega_0$@-Lagrangian~$m$@-plane~$E$ is calibrated by~$\phi_0$ 
if and only if~$\lambda(E)=1$.  For this reason, Harvey and Lawson call 
the $\phi_0$-calibrated $m$@-planes {\it special Lagrangian}.  
In particular, an $\omega_0$@-Lagrangian $m$-plane~$E$ is 
$\phi_0$@-calibrated with respect to one of its two possible orientations 
if and only if $\iota_E^*\bigl(\psi_0\bigr) = 0$, 
where $\psi_0=\Im\Upsilon_0$.  

\subsubhead 0.2.2. $\SU(m)$@-structures \endsubsubhead
Let~$\pi:P\to M$ be an~$\SU(m)$-structure on a manifold~$M$ 
of dimension~$2m$.  The elements of~$P_x=\pi^{-1}(x)$ are 
isomorphisms~$u:T_xM\to\bbC^m$ and~$\pi:P\to M$ is principal 
right $\SU(m)$@-bundle over~$M$ with right action given 
by~$u\cdot a = a^{-1}\circ u$ for~$a\in\SU(m)$.  
Then~$P$ defines (and, indeed, is defined by) the metric~$g$, 
the 2-form~$\omega$, and the complex $n$-form~$\Upsilon$ defined by
$$
\aligned
g_x &= u^*(g_0)\\
\omega_x &= u^*\bigl(\omega_0\bigr)\\
\Upsilon_x &= u^*\bigl(\Upsilon_0\bigr)\\
\endaligned
\qquad \text{for any~$u\in P_x$ and~$x\in M$.}
$$

By the above results of Harvey and Lawson, the 
$m$@-form~$\phi = \Re\Upsilon$ is a calibration if it is closed. 
Moreover, an $m$@-plane~$E\subset T_xM$ is $\phi$@-calibrated if and only
if there exists a~$u\in P_x$ so that~$u(E) = \bbR^m\subset\bbC^m$.
Such an~$E$ is an $\omega$@-Lagrangian $m$@-plane and an $\omega$@-Lagrangian 
$m$-plane~$E$ is $\phi$@-calibrated with respect to one of its two possible 
orientations if and only if $\iota_E^*\bigl(\psi\bigr) = 0$, 
where $\psi = \Im\Upsilon$.

\subsubhead 0.2.3. Special Lagrangian submanifolds \endsubsubhead
The $\phi$@-calibrated submanifolds~$L^m\subset M^{2m}$
are said to be {\it special Lagrangian}. In string theory,
especially when~$m=3$ (for example, see~\cite{SYZ}), 
these submanifolds are also known as {\it BPS} or {\it supersymmetric\/} 
cycles.  Their geometry plays an important role in the understanding
of mirror symmetry.  However, I will not attempt a discussion of these
applications in this article.

When~$m=2$, a special Lagrangian~$N^2\subset M^4$ is just
an almost-complex curve for the almost-complex structure
on~$M$ whose complex volume form is~$\omega + i\psi$.

When~$m>2$, the generic $\SU(m)$@-structure on~$M^{2m}$ will
not admit any special Lagrangian submanifolds, even locally.  However,
if the ideal~$\cI\subset\Omega^*(M)$ generated algebraically 
by~$\omega$ and~$\psi$ is differentially closed, 
and, moreover, this ideal is real analytic,
then the Cartan-K\"ahler theorem can be invoked~\cite{HL} to show that any
real analytic submanifold~$W^{m-1}\subset M$ that is an integral
manifold of~$\omega$ lies in a special Lagrangian submanifold~$X^m\subset M$.  

The case of most interest is that of an~$\SU(m)$-structure with
$d\omega = d\Upsilon = 0$.  The structure~$(M,g,\omega)$ is then a 
K\"ahler structure and~$\Upsilon$ is a holomorphic volume form
that is $g$-parallel.   In this case, 
the structure~$P$ is said to be {\it Calabi-Yau}.%
\footnote{So-named because of Yau's solution~\cite{Ya} 
of the Calabi conjecture.} By abuse of language, the 
data~$(M,\omega,\Upsilon)$ are said to constitute a Calabi-Yau manifold.
For a Calabi-Yau manifold, the underlying metric~$g$ is Ricci-flat,
so that results of DeTurck and Kazdan~\cite{DK}, imply that the metric 
and the forms~$\omega$ and~$\Upsilon$ are real analytic with respect to
the real analytic structure on~$M$ induced by its complex structure.
Moreover, since a special Lagrangian submanifold~$L^m\subset M^{2m}$
is minimal, it is necessarily real analytic, as long as it is~$C^1$.

\subsubhead 0.2.4. Moduli spaces \endsubsubhead
When~$L^m\subset M^{2m}$ is a compact special Lagrangian submanifold, 
McLean~\cite{Mc} showed that the moduli space of nearby 
special Lagrangian submanifolds has a simple description:

\proclaim{Theorem} {\smc (McLean)} Let~$\pi:P\to M^{2m}$ be a Calabi-Yau
structure on~$M$.  For any closed special Lagrangian submanifold~
$L^m\subset M^{2m}$, the $L$@-component of the moduli space of 
closed special Lagrangian submanifolds of~$M$ is a smooth manifold of 
dimension~$b_1(L)$.
\endproclaim

Here is an outline of McLean's argument.  Since~$L$ is a
Lagrangian submanifold of the symplectic manifold~$(M,\omega)$,
the Darboux-Weinstein theorem implies that there exists an $\epsilon>0$
and a symplectomorphism~$e:B_\epsilon(T^*L)\to N$ between the 
$\epsilon$@-disk bundle~$B_\epsilon(T^*L)$ in the cotangent bundle to~$L$ 
(using the induced metric on~$L$) and a tubular neighborhood~$N$ of~$L$ 
in~$M$ in such a way that~$e(0_p) = p$ for all~$p\in L$.  Then
any Lagrangian immersion~$a:L\to M$ that is sufficiently $C^1$@-close 
to~$\iota:L\to M$ is of the form~$a = e\circ\alpha$ where~$\alpha$ is
a small closed 1-form on~$L$, regarded as a section of~$T^*L$.
The mapping~$D$ from small $C^{k,\delta}$-sections of~$T^*L$ 
to~$C^{k-1,\delta}$~$m$-forms on~$L$
defined by
$$
D\alpha = {(e\circ\alpha)}^*(\psi) = \alpha^*\bigl(e^*\psi\bigr)
$$
is a first order nonlinear differential operator that takes values 
in the exact $m$-forms on~$L$ since the $m$-form~$e^*\psi$ is closed, 
vanishes along the zero section of~$T^*L$ (since, by hypothesis,
$L$ is special Lagrangian), and $\alpha$ is homotopic to the zero section.  

By choosing~$e$ carefully, one can arrange that the linearization~$D'$ 
of~$D$ at~$\alpha=0$ is given by~$D'\alpha = d(\ast\alpha)$.  
A Banach space implicit function theorem argument 
now shows that the intersection of~$D^{-1}(0)$ with the space of
$C^{k,\delta}$ closed $1$-forms on~$L$ is a smooth, finite dimensional 
manifold whose tangent space at~$\alpha=0$ is the space of 
harmonic $1$-forms on~$L$.

Each special Lagrangian manifold~$L\subset M$ inherits an orientation
and a real analytic Riemannian metric from its immersion into~$M$.
Since the stabilizer in~$\SU(m)$ of the special Lagrangian plane~$\bbR^m
\subset\bbC^m$ is~$\SO(m)$, there is no evident further structure that
is induced on~$L$ by its inclusion into~$M$.

\subsubhead 0.2.5. Real structures and slices \endsubsubhead
Nearly all of the known explicit examples of closed special Lagrangian
manifolds arise as the fixed points of a special type of involution
of a Calabi-Yau manifold.  Given a Calabi-Yau manifold~$(M,\omega,\Upsilon)$,
a {\it real structure\/} is an involution~$r:M\to M$ that satisfies
$r^*\omega = -\omega$ and~$r^*\Upsilon = \overline\Upsilon$.  If the set of 
fixed points of a real structure on a Calabi-Yau manifold is non-empty, 
then it is a special Lagrangian submanifold.  This idea
has been used in~\cite{Br1}, \cite{Br3}, and~\cite{Ko} to construct
several explicit examples of special Lagrangian manifolds 
via algebro-geometric methods.

\subhead 0.3. The coassociative case \endsubhead 
The third example that will be discussed in this article is the
coassociative calibration, which exists in dimension~$7$.
Let~$x^1,\ldots,x^7$ be the standard linear coordinates on~$\bbR^7$
and set~$dx^{ij} = dx^i\w dx^j$ and~$dx^{ijk}=dx^i\w dx^j\w dx^k$, etc.  
Define
$$
\phi_0 = dx^{567} - dx^5\w(dx^{12}+dx^{34}) - dx^6\w(dx^{13}+dx^{42})
           - dx^7\w(dx^{14}+dx^{23}).
$$
The subgroup of~$\GL(7,\bbR)$ that preserves~$\phi_0$ will be denoted~$\G_2$. 
It is a compact, connected, simple Lie group of dimension~$14$,
so this designation is appropriate.  In fact,~$\G_2$ 
preserves the metric and orientation for which the coframe~$dx=(dx^i)$ 
is oriented and orthonormal~\cite{Br1}.  

\subsubhead 0.3.1. The calibration \endsubsubhead
Harvey and Lawson show that the $4$-form
$$
\ast\phi_0 = dx^{1234}-dx^{67}\w (dx^{12}+dx^{34})-dx^{75}\w(dx^{13}+dx^{42})
           - dx^{56}\w(dx^{14}+dx^{23})
$$ 
is a calibration on~$\bbR^7$. They call this $4$-form the
{\it coassociative calibration\/} and the 4-manifolds that it 
calibrates are said to be {\it coassociative\/}.%
\footnote{The term has its origin in the algebra of the octonions.}

They show that a $4$-plane~$E$ is calibrated by~$\ast\phi_0$ 
if and only if~$\iota_E^*(\phi_0)=0$. 
For example, the $4$-plane~$E_0$ defined by~$dx^5=dx^6=dx^7=0$
and oriented so that~$dx^{1234}$ is a positive volume form is 
calibrated by~$\ast\phi_0$.  Furthermore, they show that~$\G_2$
acts transitively on the coassociative $4$@-planes and that
the $\G_2$@-stabilizer of a coassociative~$E\subset\bbR^7$ acts 
faithfully on~$E$ as its group of orientation preserving isometries. 

If~$E\subset\bbR^7$ is calibrated by~$\ast\phi_0$ and~$E^\perp
\simeq\bbR^3$ is its orthogonal complement, then the assignment~$v\mapsto
-\iota_E^*(v\lhk\phi)$ defines an isomorphism~$E^\perp\to\Lambda^2_+(E)$
that is manifestly equivariant with respect to the action of the
$\G_2$@-stabilizer of~$E$. This is visible in the case of~$E_0$, 
since the forms
$$
\Omega_2 = dx^{12}+dx^{34},\qquad 
\Omega_3 = dx^{13}+dx^{42},\qquad 
\Omega_4 = dx^{14}+dx^{23}
$$
are~$\sqrt2$ times an orthonormal basis of~$\Lambda^2_+(E_0)$.

\subsubhead 0.3.2. $\G_2$-structures \endsubsubhead
Consider a $\G_2$@-structure~$\pi:P\to M^7$.  
Thus, the elements of~$P_x=\pi^{-1}(x)$ are linear 
isomorphisms~$u:T_xM\to\bbR^7$ and~$\pi:P\to M$ is a principal
right~$\G_2$@-bundle over~$M$ where, for~$a\in\G_2$, the right action of~$a$
is given by~$u\cdot a = a^{-1}\circ u$.  The~$\G_2$-structure~$P$
induces (and, by~\cite{Br1}, is defined by) a unique 3-form~$\phi$ on~$M$ 
satisfying the condition that~$\phi_x = u^*(\phi_0)$ for some (and hence
any)~$u\in P_x$.  Moreover,~$M$ has a unique metric~$g$ and 
orientation form~$\ast1$ for which~$u:T_xM\to\bbR^7$ is an oriented 
isometry for all~$u\in P_x$.  The 4-form~$\ast\phi$ 
satisfies~$(\ast\phi)_x = u^*(\ast\phi_0)$ for all~$u\in P_x$.  

\subsubhead 0.3.3. Coassociative submanifolds \endsubsubhead
When $\ast\phi$ is closed, it is a calibration, 
known as the {\it coassociative calibration\/} of~$P$.
The 4-manifolds it calibrates are known as {\it coassociative
submanifolds\/}.  Since the subgroup of~$\G_2$ that stabilizes a 
coassociative 4-plane~$E$ is isomorphic to~$\SO(4)$ and is represented 
faithfully on~$E$ as~$\SO(E)$, a coassociative submanifold~$N^4\subset M$ 
inherits an orientation and a Riemannian metric from~$M$, 
but no finer structure of first order.

Just as in the flat case, the coassociative submanifolds are the 
4-dimensional integral manifolds of~$\phi$ and, hence, of~$d\phi$ as well. 
Now, the generic $\G_2$@-structure that satisfies~$d(\ast\phi)=0$
will not have any coassociative submanifolds because there will
not be any 4-dimensional integral manifolds of the ideal generated
by~$\phi$~and~$d\phi$.  

However, in the special case 
in which~$d\phi=0$, the situation is different.
By a theorem of Fernandez and Gray~\cite{FG}, the form~$\phi$ is
parallel with respect to the Levi-Civita connection of~$g$ 
if and only if~$d\phi = d(\ast\phi) = 0$.  In particular, in
this case, the holonomy of~$g$ preserves~$\phi$ and so is isomorphic
to a subgroup of~$\G_2$.  By a theorem of Bonan~\cite{Be,10.64}, 
this implies that the metric is Ricci flat.  Hence, by Deturck
and Kazdan~\cite{DK}, the metric~$g$ is real analytic in
harmonic coordinates.  Since~$\phi$ is harmonic with respect to~$g$,
it, too, is real analytic in these coordinates.

By the same argument used by Harvey and Lawson in the flat case,
the Cartan-K\"ahler Theorem can now
be applied to show that any real analytic 3-dimensional 
integral manifold of~$\phi$ lies in a unique 4-dimensional integral 
manifold of~$\phi$. (Any $C^1$ coassociative 4-manifold~$N$ will
be real analytic anyway, since it is minimal.) 

Thus, when~$d\phi=d(\ast\phi)=0$, the coassociative submanifolds 
of~$M$ are plentiful.  When these conditions are satisfied, the
data~$(M,\phi)$ are said to constitute a $\G_2$@-manifold.

\subsubhead 0.3.4. Coassociative moduli \endsubsubhead
When~$N^4\subset M^7$ is a closed coassociative submanifold, McLean~\cite{Mc}
showed that the moduli space of nearby closed coassociative submanifolds has
a simple description:

\proclaim{Theorem} {\smc (McLean)} Let~$(M,\phi)$ be a~$\G_2$@-manifold.  
For any closed coassociative submanifold~$L^4\subset M^7$, 
the $L$@-component of the moduli space of 
closed coassociative submanifolds of~$M$ 
is a smooth manifold of dimension~$b_2^+(L)$.
\endproclaim

Here is an outline of McLean's argument.  The normal bundle~$\nu_L$ 
can be identified with~$\Lambda^2_+(TL)$ in a natural way.  Namely, 
to every normal vector field~$u$ along~$L$, one associates 
the $2$-form~$u^\flat = -\iota_L^*(u\lhk\phi)$.  From the form of~$\phi_0$
and the fact that~$\G_2$ acts transitively on the coassociative planes,
it follows that the mapping~$u\mapsto u^\flat$ defines an isomorphism 
from~$\nu_L$ to~$\Lambda^2_+(TL)$.  Let~$\cdot^\sharp:\Lambda^2_+(TL)\to\nu_L$
denote the inverse of~$\cdot^\flat:\nu_L\to\Lambda^2_+(TL)$.
Since~$L$ is compact, each sufficiently small section~$\beta$ 
of~$\Lambda^2_+(TL)$ defines a deformation~$\widehat\beta:L\to M$ 
of~$\widehat0=\iota_L:L\hookrightarrow M$ by the formula
$$
\widehat\beta(x) = \exp_x\bigl(\beta(x)^\sharp\bigr)
$$
for~$x\in L$.  Moreover, every~$f:L\to M$ that is sufficiently $C^1$@-close
to~$\iota_L$ is of the form~$f = \widehat\beta\circ\varphi$ for some
small~$\beta\in\Omega^2_+(L)$ and 
diffeomorphism~$\varphi:L\to L$.

Now, the mapping~$D$ from small sections of~$\Lambda^2_+(TL)$ to~$3$-forms 
on~$L$ defined by
$$
D\beta = {\widehat\beta}^*(\phi)
$$
is a first-order, nonlinear differential operator.  It takes values 
in the exact $3$-forms on~$L$ since the 3-form~$\phi$ is closed, 
$\iota_L^*\phi=0$, and $\widehat\beta$ is homotopic to~$\iota_L 
= \widehat0$.  The linearization of~$D$ at~$\beta=0$ is computed to be~
$d:\Omega^2_+(L)\to\Omega^3(L)$.  A Banach space implicit function theorem 
argument now shows that~$D^{-1}(0)$ is a smooth, finite dimensional manifold 
whose tangent space at~$\beta=0$ is the space of closed, self-dual 
$2$-forms on~$L$.

\subsubhead 0.3.5. Anti@-$\G_2$ involutions and slices \endsubsubhead
If~$(M,\phi)$ is a $\G_2$@-manifold, an {\it anti@-$\G_2$\/} mapping
is a map~$r:M\to M$ that satisfies~$r^*\phi = -\phi$.  If such an~$r$
is an involution, then its fixed point set consists of a collection of 
isolated points and ($4$@-dimensional) coassociative submanifolds.  Nearly
all of the explicitly known closed coassociative submanifolds are
of this kind.  There are not many explicitly known compact examples 
beyond the ones in~\cite{BSa}.  However, many of the $\G_2$@-manifolds 
proved to exist by Joyce~\cite{Jo} admit anti@-$\G_2$ involutions 
whose fixed point locus is coassociative.

\subhead 0.4. Cartan-K\"ahler theory \endsubhead
Since the Cartan-K\"ahler Theorem plays such an important role
in this article and since many readers may not be familiar with
it, I will now recall the rudiments of 
Cartan-K\"ahler theory.  For details and proofs,
the reader may consult~\cite{BCG, Chapter 3}.

Let~$M$ be a manifold and let~$\cI\subset\Omega^*(M)$
be a graded ideal in the ring~$\Omega^*(M)$ of differential forms on~$M$.
Let~$\cI^p = \cI\cap\Omega^p(M)$.

An {\it integral manifold\/} of~$\cI$ is 
a submanifold~$\iota:N\hookrightarrow M$
so that~$\iota^*\psi = 0$ for all~$\psi\in\cI$.  An {\it integral element\/}
of~$\cI$ is a $p$-plane~$E\subset T_xM$ that satisfies~$\iota_E^*\psi=0$
for all~$\psi\in\cI$.  Let~$V_p(\cI)\subset\Gr_p(TM)$ denote the (closed) 
subset consisting of the $p$@-dimensional integral elements of~$\cI$.

The tangent spaces of an integral manifold
of~$\cI$ are evidently integral elements of~$\cI$.  The fundamental goal 
of Cartan-K\"ahler theory is to find conditions under which
a given~$E\in V_p(\cI)$ can be shown to be tangent
to some $p$-dimensional integral manifold of~$\cI$.

The concepts that need to be introduced are {\it ordinary
integral element\/}, {\it polar space\/}, and {\it regular integral
element\/}.  I will take these up in turn.

Roughly speaking, an integral element~$E\in V_p(\cI)$ is ordinary if
it is a smooth point of~$V_p(\cI)$ and the $p$-forms in~$\cI$
define~$V_p(\cI)$ `cleanly' in a neighborhood of~$E$ in~$\Gr_p(TM)$.
A more precise definition will now be given.  Let~$E\in\Gr_p(T_zM)$ be
a $p$-plane. Choose coordinates~$x = (x^i)$ on a 
$z$-neighborhood~$U\subset M$ in such a way 
that~$\iota_E^*(dx^1\w\ldots\w dx^p)\not=0$.
The set~$V$ consisting of the~$p$@-planes~$F\in\Gr_p(TU)$ 
that satisfy~$\iota_F^*(dx^1\w\ldots\w dx^p)\not=0$ is
an open~$E$@-neighborhood.  
Each~$\varphi\in\Omega^p(M)$ defines a smooth 
function~$\widehat\varphi$ on~$V$ by the rule
$$
\iota_F^*(\varphi) = \widehat\varphi(F)\,\,\iota_F^*(dx^1\w\ldots\w dx^p).
$$
Since~$\cI$ is an ideal,~$V\cap V_p(\cI) = \{\, F\in V\ \vrule\  
\widehat\varphi(F) = 0,\ \forall\varphi\in\cI^p\,\}$.  

An~$E\in V_p(\cI)$ is said to be {\it ordinary\/}
if there are~$\varphi_1,\ldots,\varphi_q\in\cI^p$ 
and an open $E$@-neighborhood~$W\subset V$
such that the functions
${\widehat\varphi}_1\,,\,\ldots,\,{\widehat\varphi}_q$
have linearly independent differentials on~$W$ and, moreover,
$$
W\cap V_p(\cI) 
    = \{\ F\in W\ \vrule\ \   
         {\widehat\varphi}_1(F) =\cdots ={\widehat\varphi}_q(F) = 0\ \}.
$$
(The condition of being ordinary does not depend on the choice of the 
coordinate system~$x$.)  Let~$V_p^o(\cI)\subset V_p(\cI)$
denote the space of ordinary integral elements of~$\cI$.  
Note that $V_p^o(\cI)$ is open in~$V_p(\cI)$ (though it might be empty).

The {\it polar space\/}~$H(E)$ 
of~$E\in V_p(\cI)$ is the union of all the integral elements of~$\cI$
that contain~$E$.  If~$e_1,\ldots, e_p$ is any basis of~$E\subset T_zM$, 
then
$$
H(E) = \{\ v\in T_zM\ \vrule\ \kappa(e_1,\ldots,e_p,v)=0\  
            \forall \kappa\in\cI^{p+1}\ \}.
$$
Consequently,~$H(E)$ is a linear subspace of~$T_zM$ that contains~$E$.
Moreover, the $(p{+}1)$@-dimensional integral elements of~$\cI$ that
contain~$E$ are the $(p{+}1)$@-dimensional subspaces of~$H(E)$ 
that contain~$E$.  If this set is non-empty, it is in one-to-one 
correspondence with the points of the projective 
space~$\bbP\bigl(H(E)/E\bigr)$.  The number~$r(E) = \dim H(E){-}p{-}1$ 
is the dimension of this space of `integral enlargements' of~$E$ and will be 
referred to as the {\it extension rank\/} of~$E$.  
Unless~$E$ is a maximal integral element,~$r(E)\ge0$ .

An ordinary integral element~$E\in V_p^o(\cI)$ is said to be {\it regular\/} 
if the extension rank function~$r:V_p(\cI)\to\{-1,0,1,\ldots\}$ is locally
constant near~$E$.  An integral manifold is said to be {\it regular\/}
if all of its tangent spaces are regular integral elements.  For
a connected regular integral manifold~$N\subset M$, the extension rank
of~$N$ is defined to be the extension rank of (any one of) its tangent spaces.

I can now state the Cartan-K\"ahler theorem~\cite{BCG}.

\proclaim{Theorem} {\smc (Cartan-K\"ahler)}
Let~$M$ be a real analytic manifold and suppose that~$\cI\subset\Omega(M)$ 
is a real analytic ideal that is closed under exterior differentiation.
Let~$X^p\subset M$ be a connected, real analytic, 
regular integral manifold of~$\cI$ and suppose that its extension
rank~$r$ is nonnegative.  Let~$Z\subset M$ be a real
analytic manifold of codimension~$r$ in~$M$ that contains~$X$ 
and satisfies%
\footnote{This is just the condition that~$T_xZ$ and~$H(T_xX)$
be transverse.}
$\dim\bigl(T_xZ\cap H(T_xX)\bigr)=p{+}1$ 
for all~$x\in X$.

Then there exists a real analytic $(p{+}1)$-dimensional integral
manifold~$Y$ of~$\cI$ that satisfies~$X\subset Y\subset Z$.  Moreover,
$Y$ is locally unique in the sense that, for any real analytic 
$(p{+}1)$@-dimensional integral manifold~$Y'$ of~$\cI$ satisfying~$X\subset
Y'\subset Z$, the intersection~$Y\cap Y'$ is also a $(p{+}1)$@-dimensional 
integral manifold of~$\cI$.
\endproclaim

Verifying the regularity of an integral manifold can be arduous.  
However, there are a few criteria for regularity that simplify this task.

One very useful fact is that, when~$E\in V_p(\cI)$ is regular,
every~$E^+\in V_{p+1}(\cI)$ that contains~$E$ is ordinary.

Building on this idea, Cartan devised a test for regularity that
is usually not difficult to check in practice.  It will now be
described.  For simplicity, assume that the ideal~$\cI$ is generated 
in positive degrees, i.e., that~$\cI^0 = (0)$.
(This will be true for the ideals that appear in this article.)

An {\it integral flag\/} of length~$n$ at~$x\in M$ is 
a sequence~$F=(E_0, E_1, E_2, \ldots, E_n)$ of integral elements 
based at~$x$ that satisfy~$\dim E_i = i$ and
$$
(0)_x=E_0\subset E_1\subset\cdots\subset E_{n-1}\subset E_n\subseteq T_xM.
$$
If~$E_i$ is regular for~$i<n$, the flag~$F$ is said to be {\it regular}.
Let~$c_i = c_i(F)$ be the codimension of~$H(E_i)$ in~$T_xM$ and define
$C = C(F) = c_0  + \cdots + c_{n-1}$.

\proclaim{Theorem} {\smc (Cartan's Test)} 
Let~$F$ be an integral flag of length~$n$. Then near~$E_n$,
the space~$V_n(\cI)$ lies in a codimension~$C$ 
submanifold of~$\Gr_n(TM)$.  
Moreover,~$F$ is regular if and only if $V_n(\cI)$ 
is a smooth submanifold of codimension~$C$ in~$\Gr_n(TM)$ 
in some neighborhood of~$E_n$.
\endproclaim

\subhead 0.5. $G$@-structures and ideals \endsubhead
This section contains an account of the relation between
torsion-free $G$@-structures and differential ideals.  For details
and proofs, see~\cite{Br2}.
When indices are needed in this subsection, I will adopt the Einstein 
summation convention and will let lower case latin indices 
range from~$1$~to~$n$.

\subsubhead 0.5.1. The coframe bundle \endsubsubhead
Let~$M^n$ be a smooth $n$@-manifold and let~$\pi:F\to M$ denote
its bundle of $\bbR^n$@-valued coframes, 
i.e., an element of~$F_x = \pi^{-1}(x)$
is a vector space isomorphism~$u:T_xM\to\bbR^n$.  This is a 
principal right $\GL(n,\bbR)$@-bundle over~$M$ where the right
action is given by~$u\cdot a = a^{-1}\circ u$ for~$a\in\GL(n,\bbR)$.
For any subspace~$\eus\subseteq\eugl(n,\bbR)=M_n(\bbR)$ 
and any~$u\in F$, let $\eus_u\subset T_uF$ denote the subspace 
of~$\ker \pi'(u)$ that corresponds to~$\eus$ under the natural 
identification of~$\eugl(n,\bbR)$ with~$\ker\pi'(u)$ generated
by this right action.

\subsubhead 0.5.2. $G$@-structures \endsubsubhead
Let~$G\subset\SO(n)$ be a connected, proper, closed Lie subgroup%
\footnote{The restriction to subgroups of~$\SO(n)$ is not 
essential, but it simplifies the presentation and is all that will 
be needed in this article.} 
with Lie algebra~$\eug\subset\euso(n)$.
The quotient space~$S = F/G$ carries the structure
of a smooth manifold.  The quotient mapping~$\tau:F\to S$ and the
induced mapping~$\bar\pi:S\to M$ are both smooth fiber bundles. 

A {\it $G$@-structure\/} on~$M$ is a principal $G$@-subbundle~$P\subset F$. 
The $G$@-structures on~$M$ are in one-to-one correspondence with the
sections of~$\tau:S\to M$.  In fact, if~$\sigma:M\to S$ is a section
of~$\tau$, then~$P_\sigma = \tau^{-1}\bigl(\sigma(M)\bigr)$ is 
a~$G$@-structure on~$M$.  Conversely, every $G$@-structure on~$M$ 
is~$P_\sigma$ for some unique $\tau$@-section~$\sigma:M\to S$.

Since~$G\subset\SO(n)$ by hypothesis, every~$G$@-structure~$P=P_\sigma$ 
has an underlying Riemannian metric~$g=g_\sigma$ and orientation that is 
defined by the condition that $u:T_xM\to\bbR^n$ be an oriented isometry 
for all~$u\in P_x$ and all~$x\in M$.  The set~$P\cdot\SO(n)\subset F$ is the 
oriented orthonormal coframe bundle of the underlying metric and orientation.

\subsubhead 0.5.3. Torsion-free and flat \endsubsubhead
A $G$@-structure~$P$ on~$M$ and the section~$\sigma$ 
so that~$P = P_\sigma$ are said to be {\it torsion-free\/} if~$P$
is parallel with respect to the Levi-Civita connection of the 
underlying Riemannian metric. The condition of being torsion-free 
is~$np$ first order quasilinear PDE for the section~$\sigma:M\to S$, 
where~$p={n\choose2}-\dim G$ is the codimension of~$G$ in~$\SO(n)$.

A torsion-free $G$@-structure~$P$ and its corresponding section~$\sigma$
are said to be~{\it flat\/} if the underlying metric~$g$ is flat.  
All flat $G$@-structures are locally equivalent to the 
translation-invariant $G$@-structure~$P_0$ on~$\bbR^n$ containing
the identity coframe~$dx:T_{\bold 0}\bbR^n\to \bbR^n$.

It is not difficult to show that $\sigma:M\to S$ 
is torsion-free if and only if it is flat to first order at every point, 
i.e., every~$x\in M$ has an open neighborhood~$U$ on which there is a flat
$G$@-structure with section~$\sigma_0:U\to S$ so that~$\sigma_0(x)=\sigma(x)$
and so that the graphs~$\sigma_0(U)$ and~$\sigma(M)$ are tangent 
at~$\sigma(x)$. For this reason, torsion-free $G$@-structures are 
sometimes referred to as being {\it $1$@-flat\/}.

As examples: When~$n=2m$ and~$G = \Un(m)$, a $\Un(m)$@-structure~$\pi:P\to M$ 
is torsion-free if and only if it is a K\"ahler structure on~$M$.  When~$n=2m$
and~$G=\SU(m)$, a~$\SU(m)$@-structure~$\pi:P\to M$ is torsion-free if and only
if it is Calabi-Yau.

\subsubhead 0.5.4. Admissible groups \endsubsubhead
Let~$\Lambda^*(\bbR^n)^G\subset \Lambda^*(\bbR^n)$ be the ring 
of~$G$@-invariant, constant coefficient differential forms on~$\bbR^n$.
Say that~$G$ is {\it admissible\/} if~$G$ is the subgroup of~$\GL(n,\bbR)$
that leaves invariant all of the forms in~$\Lambda^*(\bbR^n)^G$.  

For example, $\Un(m)\subset\SO(2m)$ is not admissible since, in this case,~
$\Lambda^*(\bbR^{2m})^G$ is generated by the $2$-form~$\omega_0$ and
yet the stabilizer of~$\omega_0$ is~$\Symp(m,\bbR)\subset\GL(2m,\bbR)$, 
which properly contains~$\Un(m)$.  

However,~$\SU(m)$ is admissible, 
since in this case, $\Lambda^*(\bbR^{2m})^G$ is generated by~$\omega_0$ and 
the real and imaginary parts of~$\Upsilon_0$ 
and~$\SU(m)=\Symp(m,\bbR)\cap\SL(m,\bbC)$. According to~\cite{Br2, 
Proposition~1}, $G = \G_2\subset\SO(7)$ is also admissible. 
In this case,~$\Lambda^*(\bbR^{7})^G$ generated 
by~$\phi_0\in\Lambda^3(\bbR^{7})$ and~$\ast\phi_0\in\Lambda^4(\bbR^{7})$.

\subsubhead 0.5.5. An ideal \endsubsubhead
Any $\alpha\in\Lambda^p(\bbR^n)$ defines a $\pi$-semibasic 
$p$-form~$\hat\alpha$ on~$F$ by
$$
\hat\alpha_u(v_1,\ldots,v_p) = \alpha\bigl(\eta(v_1),\ldots,\eta(v_p)\bigr)
$$
for~$v_1,\ldots,v_p\in T_uF$.  If~$\alpha$ is~$G$@-invariant, 
then $\hat\alpha$ is invariant under the right action of~$G$ on~$F$ 
and so descends to a well-defined $p$@-form (also denoted~$\hat\alpha$) 
on~$S$.  Any section~$\sigma:M\to S$ then induces a corresponding 
$p$@-form~$\alpha_\sigma=\sigma^*\bigl(\hat\alpha\bigr)$ on~$M$.  
When~$\sigma$ is torsion-free, the form~$\alpha_\sigma$ is parallel 
with respect to the underlying Levi-Civita connection and so must be closed.  
Consequently,~$\sigma^*\bigl(d\hat\alpha\bigr)=0$.

Let~$\cI$ denote the ideal on either%
\footnote{I will rely on context to make clear which is meant
in any given situation. In fact, when~$G$ is admissible, 
the fibers of~$\tau:F\to S$ are the Cauchy 
leaves of~$\cI$~\cite{BCG, Chapter~2, Theorem~2.2}.}
$F$~or~$S$ that is generated algebraically by the closed 
forms~$d\hat\alpha$ for~$\alpha\in\Lambda^*(\bbR^n)^G$.  
By the above discussion, the graph of a torsion-free 
section~$\sigma:M\to S$ is necessarily an integral manifold of~$\cI$. 
The converse is not always true.  

For example, when~$G = \Un(m)$, the ideal~$\cI$ is generated by the 3-form
$d\widehat{\omega_0}$ and the condition that~$\sigma(M)$ be an integral
manifold of~$\cI$ is just that~$\omega_\sigma = \sigma^*(\widehat{\omega_0})$
be closed.  Such structures are sometimes referred to in the literature
as {\it almost K\"ahler\/} since the underlying almost complex structure
of such a structure need not be integrable.  

Even when $G$ is admissible,
the closure of the forms~$\alpha_\sigma$ for~$\alpha\in\Lambda^*(\bbR^{n})^G$
need not imply that~$P_\sigma$ is torsion-free.  For example,~
$G=\Symp(2)\Symp(1)\subset\SO(8)$ is admissible and the 
ring~$\Lambda^*(\bbR^{8})^G$ is generated by a single~$4$@-form~$\Phi$.
However, it can be shown that the closure of~$\Phi_\sigma$
does not imply that the~$G$@-structure~$P_\sigma$ is torsion-free.%
\footnote{ The closure of~$\Phi_\sigma$ is only ${8\choose 5} = 56$ equations
on~$\sigma$, but the torsion-free condition is $8\cdot(28-13) = 120$ 
equations.  Of course, this, by itself, is not conclusive, since, 
{\it a priori\/}, some combination of the derivatives of the $56$ 
equations~$d\Phi_\sigma=0$ could imply the remaining~$64$ equations.  
However, a Cartan-K\"ahler analysis of this system carried out jointly 
by Dominic Joyce and myself (in 1994, but so far unpublished) 
shows that this does not happen.}

However, for the examples important in this article, those 
of~$\SU(m)\subset\SO(2m)$ and~$\G_2\subset\SO(7)$, the graph of a 
section~$\sigma:M\to S$ is an integral manifold of~$\cI$ if and only
if~$\sigma$ is torsion-free.  For~$\SU(m)$, this is well-known (and,
in any case, easy to prove). 
For~$\G_2$, this is the theorem of Fernandez and Gray mentioned earlier.

For each~$k\le n$, let~$V_k(\cI,\bar\pi)\subset\Gr_p(TS)$ denote 
the space of $k$-dimensional integral elements~$E\subset T_sS$ of~$\cI$ 
that are $\bar\pi$@-transverse, i.e., 
the projection~$\bar\pi':E\to T_{\bar\pi(s)}M$ is injective.  Similarly,
let~$V_k(\cI,\pi)\subset\Gr_p(TS)$ denote the space of $k$-dimensional
integral elements~$E\subset T_uF$ of~$\cI$ that are $\pi$@-transverse. 
Evidently, the map~$\tau':TF\to TS$ induces a surjective 
mapping~$\tau':V_k(\cI,\pi)\to V_k(\cI,\bar\pi)$, with~$\tau'(E_1)=\tau'(E_2)$
for~$E_1,E_2\subset T_uF$ if and only if~$E_1+\eug_u = E_2+\eug_u$.  
Note also that
$$
H(E) = (\tau')^{-1}\bigl(H(\tau'(E)\bigr)
$$
for all~$E\in V_k(\cI,\pi)$.  In particular, $c(E) = c\bigl(\tau'(E)\bigr)$
for such integral elements. 

\subsubhead 0.5.6. Strong admissibility \endsubsubhead
Examining the form of the generators~$d\hat\alpha$ of~$\cI$, one can  
show that~$V_n(\cI,\bar\pi)$ is a submanifold of~$\Gr_n(TS)$,
that~$V_n(\cI,\pi)$ is a submanifold of~$\Gr_n(TF)$, and that
$$
\codim\bigl(V_n(\cI,\bar\pi), \Gr_n(TS)\bigr)
= \codim\bigl(V_n(\cI,\pi), \Gr_n(TF)\bigr).
$$

Now,~$V_n(\cI,\bar\pi)$ contains the set of tangent spaces to the 
graphs of local torsion-free sections of~$S$ (which is
the same as the set of tangent spaces to the local flat sections
of~$S$).  An admissible~$G\subset\SO(n)$ is said to be
{\it strongly admissible\/} if $V_n(\cI,\bar\pi)$ consists exactly
of these tangent spaces.  

Thus, if~$G$ is strongly admissible, any section~$\sigma:M\to S$ 
whose graph~$\sigma(M)\subset S$ is an integral manifold of~$\cI$ 
is torsion-free.

When~$G\subset\SO(n)$ is strongly admissible, the space~$V_n(\cI,\bar\pi)$
is a submanifold of~$\Gr_n(TS)$ of codimension~$np$ while~
$V_n(\cI,\pi)$ is a submanifold of~$\Gr_n(TF)$ with the same codimension.
When~$G$ is not strongly admissible, this codimension will be strictly
less than~$np$.

It is not difficult to show that~$\SU(m)\subset\SO(2m)$ 
and~$\G_2\subset\SO(7)$ are strongly admissible (which implies the
theorem of Fernandez and Gray). (Some other strongly admissible groups 
are~$\Symp(m)\subset\SO(4m)$ ($m\ge1$), 
~$\Symp(m)\Symp(1)\subset\SO(4m)$ ($m\ge 3$), 
and $\Spin(7)\subset\SO(8)$.  As already remarked, 
$\Symp(2)\Symp(1)\subset\SO(8)$ is not strongly admissible.)

\subsubhead 0.5.7. Canonical flags and regular presentations \endsubsubhead
An $\bbR^n$@-valued 1-form~$\eta$ is defined on~$F$
by~$\eta(v) = u\bigl(\pi'(v)\bigr)$ for all~$v\in T_uF$.
I will usually express this form in components as~$\eta = (\eta^i)$.

Any~$E\in V_n(\cI,\pi)$ at~$u\in F$ is the terminus of 
a canonical flag~$F$ defined by
$$
E_i = \{ v\in E\ \vrule\ \eta^j(v)=0\ \text{for $j>i$}\ \}.
$$
For this canonical integral flag, the sequence of polar spaces is easy to
compute:  For~$0\le k\le n$, let~$\iota_k:\bbR^k\to\bbR^n$ denote the 
natural inclusion and set
$$
\euh_k = \{\ x\in M_n(\bbR)\ \vrule\ 
                \iota_k^*(x.\alpha) = 0, 
            \forall\alpha\in\Lambda^*(\bbR^n)^G\ \} .
$$
Then computation shows that~$H(E_k) = E + (\euh_k)_u$ for~$0\le k\le n$.
By Cartan's Test
$$
\sum_{i=0}^{n-1} \codim\bigl(\euh_i,M_n(\bbR)\bigr) 
  \le \codim\bigl(V_n(\cI,\pi),\Gr_n(TF)\bigr).
$$
The group~$G$ is said to be {\it regularly presented\/} in~$\SO(n)$ 
if equality holds.

Now, from the definition,~$\euh_k$ always contains
both~$\euh_{k+1}$ and~$M_{n,k}(\bbR)$, the space
of $n$@-by@-$n$ matrices whose first~$k$ columns are all zero.
When~$G$ is admissible, $\euh_n = \eug$ and it is
not difficult to calculate that~$\euh_{n-1}=\eug+M_{n,n-1}(\bbR)$.
For~$k<n{-}1$ however, the space~$\euh_k$ is not so readily computed.
In fact, the dimensions of the 
spaces~$\euh_k$ can depend on~$\eug$ itself and not just on the conjugacy
class of~$\eug$ in~$\euso(n)$.

For example, for any integer~$m>0$, consider the $2m$@-by@-$2m$ matrix
$$
\text{J}_m = \pmatrix 0_m & \text{I}_m\\ -\text{I}_m & 0_m \endpmatrix
$$
and the $2m$@-by@-$2m$ matrix
$$
\text{J}^*_m = \pmatrix 
                \text{J}_1 & 0_2 & \cdots & 0_2\\
                0_2 & \text{J}_1 & \cdots & 0_2\\
                \vdots & \vdots & \ddots & \vdots\\
                0_2 & 0_2 & \cdots & \text{J}_1 \\
                \endpmatrix.
$$
These two matrices are conjugate in~$\Or(2m)$.  Let
$$
\openup1\jot
\align
\eusu(m)&=\{ x\in\euso(2m)\ \vrule\ 
                   x\text{J}_m=\text{J}_mx,\ \ \tr(\text{J}_mx) = 0\ \}\\
^*\!\eusu(m)&=\{ x\in\euso(2m)\ \vrule\ 
               x\text{J}^*_m=\text{J}^*_mx,\ \ \tr(\text{J}^*_mx) = 0\ \}\\
\endalign
$$
and let~$\SU(m)$ and $^*\!\SU(m)$ denote the corresponding (conjugate) 
Lie subgroups.  It can be shown that~$\SU(m)$ is regularly 
presented, but that, for~$m>2$, the group~$^*\!\SU(m)$ is not.  
It is for this reason that I  defined~$\SU(m)$ in~\S0.2 
to be the former group and not the latter.

\subsubhead 0.5.8. An existence result \endsubsubhead
If~$G$ is regularly presented, then by Cartan's Test 
every~$E\in V_n(\cI,\pi)$ is the terminus of a regular flag, namely~$F$.  
Moreover, setting~$\bar E = \tau'(E)$ and considering the flag~$\bar F$ 
defined by~$\bar E_i = \tau'(E_i)$, it follows that the 
flag~$\bar F$ is also regular.  

The smooth manifold~$M$ has an underlying real analytic structure.%
\footnote{This real analytic structure is essentially unique by a 
theorem of Grauert~\cite{Gr}.}
The bundles~$P$ and~$S$ and the ideal~$\cI$ then inherit real 
analytic structures from that of~$M$.  Since~$\cI$ is differentially
closed by construction, the Cartan-K\"ahler theorem has the following
consequence:

\proclaim{Corollary}
If~$G\subset\SO(n)$ is conjugate to a regularly presented subgroup 
of~$\SO(n)$, then every~$E\in V_n(\cI,\bar\pi)$ is tangent to the graph 
of some local section~$\sigma:U\to S$ defined on some 
open~$U\subset M$ and satisfying~$\sigma^*(\cI) = 0$. \qed
\endproclaim

When~$G$ is also strongly admissible, the section~$\sigma$ will be
torsion-free.%
\footnote{However, this result is of interest even when~$G$ is not
strongly admissible.  For example, when~$G = \Symp(2)\Symp(1)\subset\SO(8)$,
it can be shown that~$G$ is conjugate to a regularly presented subgroup.
This is the heart of the calculation by Joyce and myself referred to 
in an earlier footnote.}

\head 1. Special Lagrangian Realization  \endhead

\proclaim{Theorem 1}
Let~$(L^3,g)$ be a closed, real analytic, oriented Riemannian manifold.  
Then there exists a Calabi-Yau $3$-fold~$(N^6,\omega,\Upsilon)$ and an
isometric embedding~$\iota:L\to N$ whose image is a special Lagrangian
$3$-manifold in~$N$.  Moreover, $(N^6,\omega,\Upsilon)$ and $\iota$
can be chosen so that~$\iota(L)$ is the fixed
locus of a real structure~$r:N\to N$.
\endproclaim

\demo{Proof} First, I will examine the differential system~$\cI$
constructed in \S0.5 for the group~$G = \SU(3)$ and show that
it is regularly presented.  The summation convention is still in force
and lower case Latin indices now range from 1 to 6.

As in \S0.2, let $\SU(3)\subset\SO(6)$ be the connected subgroup
of dimension~$8$ whose Lie
algebra consists of the matrices of the form
$$
\pmatrix a & b\\ -b & a \endpmatrix
$$
where~$a$ is a skewsymmetric 3-by-3 matrix and~$b$ is a traceless
symmetric 3-by-3 matrix.  This~$\SU(3)$ can also be defined as the 
simultaneous stabilizer in~$\GL(6,\bbR)$ of the forms $\omega_0$
and~$\Upsilon_0$ in~$\Lambda^*(\bbR^6)$ defined by the formulae
$$
\align
\omega_0 &= dx^{14} + dx^{25} + dx^{36}\,,\\
\Upsilon_0 &= (dx^1+\imath\,dx^4)\w 
            (dx^2+\imath\,dx^5)\w(dx^3+\imath\,dx^6)\\
&= dx^{123} - dx^{156} + dx^{246}
    - dx^{345} + \imath\,\bigl(dx^{246} + dx^{234} 
          - dx^{135} - dx^{456}\bigr).\\
\endalign
$$
Consequently,~$G=\SU(3)$ is admissible. 

Since~$\Lambda^*(\bbR^6)^G$ has no forms of degree~$1$, 
$\euh_0=\euh_1 = M_6(\bbR)\simeq\bbR^{36}$.  
Moreover,
$$
\iota_2^*\bigl(x.\omega_0\bigr) = (x^4_2 - x^5_1)\,dx^1\w dx^2,
$$
so~$\euh_2 = \{x\in M_6(\bbR)\ \vrule\ x^4_2 - x^5_1 = 0\}
\simeq\bbR^{35}$.  Similarly,~$\iota_3^*\bigl(x.\omega_0\bigr) 
= \iota_3^*\bigl(x.\Upsilon_0\bigr)= 0$ if and only if
$$
x^4_2 - x^5_1 = x^5_3-x^6_2 = x^6_1-x^4_3 =
  x^1_1+x^2_2+x^3_3 = x^4_1+x^5_2+x^6_3 = 0,
$$
so~$\dim\euh_3=31$.  (Note that $\euh_3$ is defined
by at most $5 = {3\choose2}+2\cdot{3\choose3}$ linearly independent 
equations anyway.)  Further, $\euh_4\subset\euh_3$ is defined by the
additional nine equations
$$
\align
x^1_1+x^1_1=x^2_1+x^1_2=x^3_1+x^1_3&=0,\\
x^1_4+x^4_1=x^2_4+x^5_1=x^3_4+x^6_1&=0,\\
x^4_4-x^1_1=x^5_4-x^2_1=x^6_4-x^3_1&=0,\\
\endalign
$$
so~$\dim\euh_4 = 22$.  As already remarked, 
the results $\dim\euh_5 = 8+6 = 14$ and~$\dim\euh_6 = 8$ 
follow from the general considerations of \S0.5.

Now let $M$ be any $6$@-manifold and let~$\pi:F\to M$ and~$S = F/\SU(3)$ 
be its coframe bundle and~$\SU(3)$@-structure bundle as in~\S0.5. 
Let~$\eta = (\eta^i)$ be the tautological $\bbR^6$@-valued 1-form on~$F$.  
The ideal~$\cI$ on~$F$ or~$S$
is generated by~$d(\widehat{\omega_0})$ 
and~$d\bigl(\widehat{\Upsilon_0}\bigr)$.

By the calculations above and the discussion in~\S0,5,
any~$E\in V_6(\cI,\pi)$ has a canonical flag~$\{E_i\,\vrule\,0\le i\le 6\}$
and the codimension of~$H(E_i)$ is~$c_i$, where
$$
(c_0,c_1,\ldots,c_6) = (0,0,1,5,14,22,28).
$$

Since~$\SU(3)$ is strongly admissible, $V_6(\cI,\pi)$ 
has codimension~$6\cdot(15-8) = 42$ in~$\Gr_6(TS)$. Since
$$
c_0 + c_1 + c_2 + c_3 + c_4 + c_5 = 42,
$$
it follows that~$\SU(3)$ is regularly presented, and the canonical flag
associated to any integral element~$E\in V_6(\cI,\pi)$ is regular.

Now I want to show how one can use the Cartan-K\"ahler Theorem to
produce the desired~$M$ of Theorem~1.

By a theorem of~Wu~\cite{MS}, an orientable $3$-manifold is
smoothly parallelizable.  Bochner~\cite{Bo} proved that on a 
closed, real analytic Riemannian manifold, the real analytic
differential forms are dense in the smooth differential forms
endowed with the uniform topology.  It follows that a real 
analytic orientable Riemannian $3$-manifold is 
real analytically parallelizable.  

Consequently, Gramm-Schmidt orthonormalization
can be invoked to construct real analytic 
1-forms~$\omega_1$, $\omega_2$, $\omega_3$ on~$L$ so that
$$
g = {\omega_1}^2+{\omega_2}^2+{\omega_3}^2
$$
and so that~$\omega_1\w\omega_2\w\omega_3$ is a positive 
volume form on~$L$ for the given orientation.

Let $M = L\times\bbR^3$ and let~$y^1$, $y^2$, and $y^3$ be linear 
coordinates on the second factor and regard them as functions
on~$M$ via pullback from its projection onto the second factor.
For notational simplicity, I will identify~$L$ with~$L\times{\bold0}$
from now on.  

The 1-forms~$(\omega_1,\omega_2,\omega_3,dy^1,dy^2,dy^3)$ 
define a real analytic parallelization of~$M$.
As before, let~$\pi:F\to M$ be the coframe bundle on~$M$
and let~$\eta$ be the canonical~$\bbR^6$@-valued 1-form.

There exists a unique real analytic map~$g:F\to\GL(6,\bbR)$ 
for which
$$
\pmatrix
   \eta^1\\\eta^2\\\eta^3\\\eta^4\\\eta^5\\\eta^6\\
\endpmatrix
= g^{-1}\,
\pmatrix
   \omega_1\\\omega_2\\\omega_3\\dy^1\\dy^2\\dy^3\\
\endpmatrix.
$$
Thus,~$(\pi,g):F\to M\times\GL(6,\bbR)$ is a real analytic trivialization
of~$F$, regarded as a~$\GL(6,\bbR)$@-bundle.  To simplify notation, 
I will usually identify~$F$ with~$M\times\GL(6,\bbR)$ 
without explicitly noting the identification.

Theorem~1 will be proved by applying the Cartan-K\"ahler theorem 
to construct a real analytic 6-dimensional integral 
manifold of~$\cI$ that projects diffeomorphically onto a neighborhood~$N$
of~$L\subset M$ so that the induced Calabi-Yau
structure on~$N$ has the properties that, first,
$L $ is a special Lagrangian submanifold of~$N$ and, second, 
the induced metric and orientation agree with the given ones on~$L$.

To prove the last statement in Theorem~1, I will construct an involution 
of~$M$ and a covering involution of~$F$.  
By abuse of notation, I will use the same letter~$r$ for each.
The formula for~$r$ acting on~$M=L\times\bbR^3$ is~$r(p,y) = (p,-y)$;
when it acts on~$F$ the formula is
$$
r(p,y,g) = (p,-y,RgR),\qquad\qquad
\text{where}\ \  R = \pmatrix \text{I}_3&0\\ 0&-\text{I}_3\endpmatrix.
$$
Now,~$R = R^{-1}\in\Or(6)$ does not belong to~$\SU(3)$, 
but it is nevertheless true that~$R\,\SU(3)\,R = \SU(3)$.
Consequently, $r$ preserves the fibers of~$\tau:F\to S$, and thus
induces an involution of~$S$, which I will continue to denote by~$r$.

Since~$r^*\eta = R\eta$ and since
$$
\align
\widehat{\omega_0} &= \eta^1\w\eta^4 + \eta^2\w\eta^5 + \eta^3\w\eta^6\,,\\
\widehat{\Upsilon_0} &= (\eta^1+\imath\,\eta^4)\w 
            (\eta^2+\imath\,\eta^5)\w(\eta^3+\imath\,\eta^6).\\
&= \eta^1\w\eta^2\w\eta^3 - \eta^1\w\eta^5\w\eta^6 + \eta^2\w\eta^4\w\eta^6
    - \eta^3\w\eta^4\w\eta^5 \\
&\qquad\qquad + \imath\,\bigl(\eta^1\w\eta^2\w\eta^6 + \eta^2\w\eta^3\w\eta^4 
          - \eta^1\w\eta^3\w\eta^5 - \eta^4\w\eta^5\w\eta^6\bigr),\\
\endalign
$$
it follows that
$$
r^*\bigl(\widehat{\omega_0}\bigr) = {}-\widehat{\omega_0}\quad
\text{and}\qquad 
r^*\bigl(\widehat{\Upsilon_0}\bigr) = \overline{\widehat{\Upsilon_0}}.
$$
In particular,~$r^*\cI = \cI$, so that $r$ takes integral
manifolds of~$\cI$ to integral manifolds of~$\cI$.

The integral manifold of~$\cI$ to be constructed in the course
of the proof will be invariant under the involution~$r$. It will
then follow that forms~$(\omega,\Upsilon)$ induced by the Calabi-Yau 
structure on the neighborhood~$N$ of~$L$ will satisfy~$r^*(\omega)=-\omega$ 
and~$r^*\Upsilon=\overline\Upsilon$. Thus~$r:N\to N$ is a real structure 
on~$(N,\omega,\Upsilon)$ and~$L\subset N$ is its fixed locus.

To begin the construction of the desired 6-dimensional integral 
manifold, define a lifting~$f_3:L\to S$ by
taking it to be of the form
$$
f_3(p) = \tau\bigl(p,{\bold0},\text{I}_6\bigr).
$$
Then~$f_3^*(\widehat{\omega_0}) = 0$ and~$f_3^*(\widehat{\Upsilon_0}) 
= \omega_1\w\omega_2\w\omega_3$.
Consequently, $f_3^*(d\widehat{\omega_0}) = f_3^*(d\widehat{\Upsilon_0}) = 0$, 
so $X_3 = f_3\bigl(L \bigr)$ is an integral manifold of~$\cI$.
Note that~$X_3$ lies in the fixed locus of~$r$.   Moreover,
the tangent spaces to~$X_3$ are the projections to~$S$
of the tangents to the lifting~$p\mapsto (p,{\bold0},\text{I}_6)\in F$,
so these are of the type of~$E_3$ in a canonical (regular) flag.  It follows
that these tangent spaces are all regular and their polar spaces have 
dimension~$29$ (the minimum possible).
Thus, the computation of~$\euh_3$ shows that~$X_3$ is a regular integral
manifold of~$\cI$ with extension rank~$25$.

Let~$W_5\subset M_6(\bbR)$ denote the $5$-dimensional subspace
consisting of the matrices of the form
$$
x = 
\pmatrix 
x_5&0&0&0&0&0\\
0&x_5&0&0&0&0\\
0&0&x_5&0&0&0\\
x_4&x_3&-x_2&0&0&0\\
-x_3&x_4&x_1&0&0&0\\
x_2&-x_1&x_4&0&0&0\\
\endpmatrix\,.
$$
Note that~$W_5\cap\euh_3=(0)$. Since~$\euh_3$ contains~$\eusu(3)$,
the affine space~$\text{I}_6+W_5$ intersects~$\SU(3)$ transversely 
at~$\text{I}_6\in\GL(6,\bbR)$.  Thus, there is a 
$\bold0$@-neighborhood~$U_5\subset W$ so that the map~$U_5\times\SU(3)\to
\GL(6,\bbR)$ defined by~$(x,a) \mapsto (\text{I}_6+x)a$ is an embedding.
Note also that~$W_5$ is invariant under conjugation by~$R$.

Define a $9$-dimensional submanifold~$Z_3\subset S$ by the rule
$$
Z_3 = \{ \tau\bigl(p,(y^1,0,0),\text{I}_6{+}x\bigr)\ 
          \vrule\ p\in L,\ y^1\in\bbR,\ x\in U_5\ \}.
$$
Note that~$X_3\subset Z_3$ and that~$Z_3$ is $r$@-invariant.  By
the description of~$\euh_3$, it follows that~$H(T_xX_3)$ and~$T_xZ_3$ meet
transversely along~$X_3$, so that $E_x = H(T_xX_3)\cap T_xZ_3$
is a $4$-dimensional integral element of~$\cI$ for all~$x\in X_3$
that is $\bar\pi$@-transverse.  
In fact, its $\bar\pi$@-projection is spanned by the corresponding 
tangent space to~$L$ plus the vector $\partial/\partial y^1$.

By the Cartan-K\"ahler theorem, there is a real analytic $4$-dimensional
integral manifold~$Y_4$ of~$\cI$ that satisfies~$X_3\subset Y_4\subset Z_3$.
Since~$X_3$ and~$Z_3$ are $r$@-invariant,~$r(Y_4)$ is also an integral 
manifold of~$\cI$ and satisfies~$X_3\subset r(Y_4)\subset Z_3$.  
By the uniqueness part of the Cartan-K\"ahler Theorem,~$X_4=Y_4\cap r(Y_4)$ 
is also a 4-dimensional integral manifold of~$\cI$ and is 
manifestly~$r$@-invariant.

A neighborhood of~$X_3$ in~$X_4$ projects diffeomorphically onto
a neighborhood~$N_4$ of~$L $ in~$L\times\bbR^1\subset M$.  
By shrinking~$X_4$ if necessary, it can be supposed that~$X_4$ 
is the graph of a section of~$S$ along such an 
open~$N_4\subset L\times\bbR^1$, so assume this.
Furthermore, for~$x\in X_3$, the integral elements~$T_xX_4$ are
of type~$E_4$ in a canonical flag and hence must be regular.
By shrinking~$X_4$ again if necessary, it can be assumed that~$X_4$ 
is a connected, regular integral manifold of~$\cI$ containing~$X_3$
and invariant under~$r$. By the calculation of~$\euh_4\simeq\bbR^{22}$, 
it follows that~$X_4$ has extension degree~$15$.  

Thus, in order to extend~$X_4$ to a $5$-dimensional integral manifold, 
one needs to find a restraining manifold~$Z_4\subset S$ 
that contains~$X_4$, is of codimension~$15$, meets the polar 
spaces of the tangent planes to~$X_4$ transversely, 
and is $r$@-invariant. This is done as follows: 

Let~$W_{14}\subset M_6(\bbR)$ denote the $14$-dimensional subspace
containing~$W_5$ and consisting of the matrices of the form
$$
x = 
\pmatrix 
x_5+x_{10}&x_{11}&x_{12}&x_8&0&0\\
-x_{14}&x_5+x_9&0&x_6&0&0\\
-x_{13}&0&x_5+x_9&x_7&0&0\\
x_4&x_3&-x_2&x_9+x_{10}&0&0\\
-x_3+x_6&x_4-x_8&x_1&x_{11}+x_{14}&0&0\\
x_2+x_7&-x_1&x_4-x_8&x_{12}+x_{13}&0&0\\
\endpmatrix\,.
$$
Note that $W_{14}\cap\euh_4=(0)$.  Since~$\euh_4$ 
contains~$\eusu(3)$, the affine space~$\text{I}_6+W_{14}$
intersects~$\SU(3)$ transversely 
at~$\text{I}_6\in\GL(6,\bbR)$.  Thus, there is a 
$\bold0$@-neighborhood~$U_{14}\subset W_{14}$ so that 
the map~$U_{14}\times \SU(3)\to\GL(6,\bbR)$ defined 
by~$(x,a) \mapsto (\text{I}_6+x)a$ is an embedding.
Note also that~$W_{14}$ is invariant under conjugation by~$R$.

Define a $19$-dimensional submanifold~$Z_4\subset S$ by the rule
$$
Z_4 = \{ \tau\bigl(p,(y^1,y^2,0),\text{I}_6{+}x\bigr)\ 
          \vrule\ p\in L,\ y^1\in\bbR,\ x\in U_{14}\ \}.
$$
Then~$Z_3\subset Z_4$ and~$Z_4$ is $r$@-invariant.  Using the
computation of~$\euh_4$, one sees that~$H(T_xX_4)$ and~$T_xZ_4$ meet
transversely along~$X_4$, so that $E_x = H(T_xX_4)\cap T_xZ_4$
is a $5$-dimensional integral element of~$\cI$ for all~$x\in X_4$
that is $\bar\pi$@-transverse.  In fact, its $\bar\pi$@-projection 
is spanned by the corresponding tangent space to~$L$ plus the vectors
$\partial/\partial y^1$ and~$\partial/\partial y^2$.

By the Cartan-K\"ahler theorem, there is a real analytic $5$-dimensional
integral manifold~$Y_5$ of~$\cI$ that satisfies~$X_4\subset Y_5\subset Z_4$.
Since~$X_4$ and~$Z_4$ are $r$@-invariant,~$r(Y_5)$ is also an integral 
manifold of~$\cI$ and satisfies~$X_4\subset r(Y_5)\subset Z_4$.  
By the uniqueness part of the Cartan-K\"ahler Theorem,~$X_5=Y_5\cap r(Y_5)$ 
is also a 5-dimensional integral manifold of~$\cI$ and is 
manifestly~$r$@-invariant.

A neighborhood of~$X_3$ in~$X_5$ projects diffeomorphically onto
a neighborhood~$N_5$ of~$L$ in~$L\times\bbR^2\subset M$.  
By shrinking~$X_5$ if necessary, it can be supposed that~$X_5$ 
is the graph of a section of~$S$ along such an 
open~$N_5\subset L\times\bbR^2$, so assume this.
Then~$X_5$ is a regular integral manifold
of~$\cI$.  By the computation of~$\euh_5$, it has extension 
degree~$6$.

Thus, in order to extend~$X_5$ to a $6$-dimensional integral manifold, 
one needs to find a restraining manifold~$Z_5\subset S$ 
that contains~$X_5$, is of codimension~$6$ and meets the polar 
spaces of~$X_5$ transversely, and is $r$@-invariant. This is done 
as follows. 

Let~$W_{22}\subset M_6(\bbR)$ denote the $22$-dimensional subspace
containing~$W_{14}$ and consisting of the matrices of the form
$$ 
\pmatrix 
   x_5+x_{10}& x_{11}&  x_{12}&         x_8&x_{15}&0\\
-x_{14}+x_{16}&x_5+x_9+x_{17}&x_{18}& x_6-x_{15}&x_{19}&0\\
      -x_{13}& -x_{21}&x_5+x_9+x_{22}&          x_7&x_{20}&0\\
   x_4-x_{19}&    x_3&   -x_2&   x_9+x_{10}+x_{22}&x_{16}&0\\
     -x_3+x_6&x_4-x_8&x_1&x_{11}+x_{14}&x_{17}+x_{22}&0\\
      x_2+x_7&-x_1+x_{20}&x_4-x_8-x_{19}&x_{12}+x_{13}&x_{18}+x_{21}&0\\
\endpmatrix\,.
$$
Then~$W_{22}\cap\euh_5=(0)$.  Since~$\euh_5$ contains~$\eusu(3)$, 
the affine space~$\text{I}_6+W_{22}$ intersects~$\SU(3)$ transversely 
at~$\text{I}_6\in\GL(6,\bbR)$.  Thus, there is a 
$\bold0$@-neighborhood~$U_{22}\subset W$ so that 
the map~$U_{22}\times \SU(3)\to\GL(6,\bbR)$ defined 
by~$(x,a) \mapsto (\text{I}_6+x)a$ is an embedding.
Note also that~$W_{22}$ is invariant under conjugation by~$R$.

Define a $28$-dimensional submanifold~$Z_5\subset S$ by the rule
$$
Z_5 = \{ \tau\bigl(p,(y^1,y^2,y^3),\text{I}_6{+}x\bigr)\ 
          \vrule\ p\in L,\ y^1\in\bbR,\ x\in U_{22}\ \}.
$$
Then~$Z_4\subset Z_5$ and~$Z_5$ is $r$@-invariant.  
By the $\euh_5$ computation, one sees that~$H(T_xX_5)$ and~$T_xZ_5$ meet
transversely along~$X_5$, so that $E_x = H(T_xX_5)\cap T_xZ_5$
is an integral element of~$\cI$ for all~$x\in X_5$
that is $\bar\pi$@-transverse.  

By the Cartan-K\"ahler Theorem, there is a real analytic $6$-dimensional
integral manifold~$Y_6$ of~$\cI$ that satisfies~$X_5\subset Y_6\subset Z_5$.
Since~$X_5$ and~$Z_5$ are $r$@-invariant,~$r(Y_6)$ is also an integral 
manifold of~$\cI$ and satisfies~$X_5\subset r(Y_6)\subset Z_5$.  
By the uniqueness part of the Cartan-K\"ahler Theorem,~$X_6=Y_6\cap r(Y_6)$ 
is also a 6-dimensional integral manifold of~$\cI$ and is 
manifestly~$r$@-invariant.

A neighborhood of~$X_3$ in~$X_6$ projects diffeomorphically onto
an open neighborhood~$N=N_6$ of~$L$ in~$M=L\times\bbR^3$.  
Thus,~$X_6$ is the graph of some~$\sigma:N\to S$.

As has already been remarked, the section~$\sigma$ induces
a Calabi-Yau structure on~$N$ that satisfies~$r^*\omega = -\omega$
and~$r^*\Upsilon=\overline\Upsilon$. By the way~$X_3$ was chosen,
the equations
$$
\align
\omega &= \omega_1\w dy^1 + \omega_2\w dy^2 + \omega_3\w dy^3\\
\Upsilon &= 
 (\omega_1+\imath\,dy^1)\w(\omega_2+\imath\,dy^2)\w(\omega_3+\imath\,dy^3)\\
\endalign
$$
hold along the locus~$y=0$.  The underlying metric of the Calabi-Yau
structure agrees with
$$
{\omega_1}^2+{\omega_2}^2+{\omega_3}^2+(dy^1)^2+(dy^2)^2+(dy^3)^2
$$
along~$y=0$.  In particular, this metric induces the original metric~$g$
on~$L$, which, being the fixed locus of the real structure~$r$, is 
necessarily special Lagrangian.\qed
\enddemo

To conclude this section, here are a few remarks about different approaches
and generalizations.

\remark{An alternative approach}
Theorem~$1$ can be given a different proof 
that is based on constructing the Calabi-Yau structure in two stages.  

First, let~$\pi:F\to M^6$ be the coframe bundle as above and
consider the quotient bundle~$\bar\pi:J = F/\SL(3,C)\to M$ 
where~$\SL(3,C)$ is the connected subgroup of~$\GL(6,\bbR)$ 
whose Lie algebra consists of the matrices of the form
$$
\pmatrix a & b\\ -b & a \endpmatrix
$$
where~$a$ and~$b$ are traceless 3-by-3 matrices.  The $3$-form
$$
\gamma
=(\eta^1+\imath\,\eta^4)\w(\eta^2+\imath\,\eta^5)\w(\eta^3+\imath\,\eta^6)
$$
is well-defined on~$J$ and~$d\gamma$ generates a differential 
ideal~$\cJ$ on~$J$.  Let~$V_p(\cJ,\bar\pi)$ denote the
space of $\bar\pi$@-transverse integral elements. 
One can then prove the following result:

\proclaim{Proposition 1}
For any~$E\in V_p(\cJ,\bar\pi)$, the dimension of~$H(E)$ is at least~$m_p$
where
$$
(m_0,m_1,m_2,m_3,m_4,m_5,m_6) = (26,26,26,24,18,12,6)
$$
and equality holds except possibly when $p=3$ or $4$.
Moreover, equality does hold for~$E$ in a dense open subset 
of~$V_p(\cJ,\bar\pi)$ and, when it does, $E$ is regular.
\endproclaim

The proof of Proposition~1 is simpler than the corresponding
result for~$G=\SU(3)$, partly because the group~$\SL(3,\bbC)$ is 
larger than~$\SU(3)$ and partly because~$\cJ$ is algebraically 
simpler than~$\cI$.

One then uses Proposition~1 to show that for~$M=L\times\bbR^3$,
with involution~$r:M\to M$ defined by~$r(p,y) = (p,-y)$ as before,
one can find a neighborhood of~$L=L\times{\bold0}$ on which there
is a section~$j$ of~$J$ whose graph is an integral manifold 
of~$\cJ$ so that the induced holomorphic volume 
form~$\Upsilon = j^*\gamma$ satisfies~$r^*\Upsilon=\overline\Upsilon$
and so that the pullback of~$\Upsilon$ to the slice~$L$
is the oriented volume form of~$g$.  Again, this proof is
simpler than the corresponding one for~$\SU(3)$ because the polar
spaces are smaller and the transversality conditions
are easier to verify.

Now suppose that~$M$ is a complex $3$-manifold endowed with
a holomorphic volume form~$\Upsilon$.  Let~$\pi_\bbC:F_\bbC\to M$
be the $\SL(3,\bbC)$@-structure whose fiber over~$x\in M$ 
consists of the complex linear, volume preserving 
isomorphisms~$u:T_xM\to\bbC^3$.  The canonical form~$\eta$
can now be regarded as a $\bbC^3$@-valued $1$-form on~$F_\bbC$.
Let
$$
\beta = \frac\imath2\,\left(\eta^1\w\overline{\eta^1}{+}
         \eta^2\w\overline{\eta^2}{+}\eta^2\w\overline{\eta^2}\right)
\qquad
\text{where}\ \ 
\eta = \pmatrix \eta^1\\ \eta^2\\ \eta^3\\ \endpmatrix.
$$
Then~$\beta$  is invariant under right action by~$\SU(3)$ and so
is well-defined on the quotient bundle~$\bar\pi_\bbC: K=F_\bbC/\SU(3)\to M$.

One can now consider the ideal~$\cK$ generated by~$d\beta$ on~$K$
and the $\bar\pi_\bbC$@-transverse integral elements~$V_p(\cK,\bar\pi_\bbC)$. 
One can then prove

\proclaim{Proposition 2}
For any~$E\in V_p(\cK,\bar\pi_\bbC)$, the dimension of~$H(E)$ 
is at least~$m_p$ where
$$
(m_0,m_1,m_2,m_3,m_4,m_5,m_6) = (14,14,13,11,8,6,6)
$$
and equality holds except possibly when $p=3$ or $4$.
Moreover, equality does hold for~$E$ in a dense open subset 
of~$V_p(\cK,\bar\pi_\bbC)$ and, when it does, $E$ is regular.
\endproclaim

Again, the proof of Proposition~2 is simpler than that Theorem~1.

Combining this result with the same sorts of arguments made in the 
proof of Theorem~$1$ (but, again, much simpler), one can now show that
starting with~$M$ a neighborhood of~$L$ already endowed
with a holomorphic volume form for which~$r$ serves as an antiholomorphic
involution, one can construct a section of~$K$ over an~$L$@-neighborhood~$N
\subset M$ whose graph is an integral manifold of~$\cK$, is invariant under 
the appropriate involution, and is such that the induced K\"ahler form~
$\omega$ together with~$\Upsilon$ defines a Calabi-Yau structure on~$N$ 
that induces the original given metric on~$L$.

While the resulting proof is locally simpler, it is no shorter than
the `collapsed' proof of Theorem~1.  It has another disadvantage in
that it does not correspond well with the proof to be offered
in the next section about coassociative extensions.
\endremark

\remark{Higher dimensions}
On the other hand, the $2$-step proof just described does have 
the advantage that it leads to an easier proof of the generalization
of Theorem~1 to 
higher dimensions: Any {\it parallelizable\/} $m$-manifold~$L^m$ 
endowed with an orientation and a real analytic Riemannian metric
can be isometrically embedded as a special Lagrangian submanifold
of a Calabi-Yau manifold.  As in Theorem~1, one can even arrange
that the special Lagrangian submanifold be the fixed locus of
an antiholomorphic involution of the Calabi-Yau structure.
\endremark

\remark{Parallelizability}
It is an interesting question as to whether the assumption of
parallelizability is really necessary.  It is certainly necessary
for the proof via Cartan-K\"ahler theory since this method 
requires that one be able to construct the extension one 
dimension at a time.  However, this does not say whether or not
there might be some other proof that works in general.  Fortunately,
since orientable $3$-manifolds are parallelizable, this difficulty
did not arise for the proof of Theorem~1.

The reader may wonder whether or not the choice of parallelization
really has any effect, so that one might imagine that one could
do the extension locally and then expect the solutions to `patch'
appropriately.  However, a careful look at how the restraining
manifolds~$Z_3$, $Z_4$, and $Z_5$ are constructed shows that
the choice of parallelization does indeed have an effect on the
solution.  I do not see any way of choosing a family of 
local parallelizations whose corresponding extensions would 
agree on overlaps.

Also, the reader may have noticed that the only use of the assumption
that~$L$ be closed was to derive the existence of a real analytic
parallelization, which is necessary for the Cartan-K\"ahler proof.  
It seems reasonable that any orientable, real analytic $3$-manifold 
would be real analytically parallelizable.   I have not been able to 
find or construct a convincing proof of this, but 
perhaps the results in~\cite{Sh} are applicable.  If so, then
Theorem~1 can be strengthened  by removing the hypothesis that~$L$
be compact.
\endremark 

\remark{Prescribing the second fundamental form}
The reader may wonder whether the method can be
used to find extensions in which the special Lagrangian manifold~$L^3$
is not totally geodesic in~$M^6$.  The answer is `yes'.  In fact,
one can prescribe the `second fundamental form' arbitrarily subject
to some obvious necessary conditions.

More precisely, the situation can be described as follows:  
For any Lagrangian submanifold~$L^m\subset M^{2m}$, the 
normal bundle to~$L$ can be canonically identified with~$T^*\!L$,
using the metric and symplectic structure.  Thus, the second
fundamental form~$\text{I\!I}$ can be regarded as a section
of~$T^*\!L\otimes S^2(T^*\!L)$.  Because~$L$ is Lagrangian,
it turns out that~$\text{I\!I}$ actually takes values in~$S^3(T^*\!L)
\subset T^*\!L\otimes S^2(T^*\!L)$.  If~$L$ is special Lagrangian,
then it is minimal, so that the trace of~$\text{I\!I}$
with respect to the induced metric~$g$ vanishes.  I.e., $\text{I\!I}$
must take values in~$S^3_0(T^*\!L)$.  

It then turns out that, given a real analytic metric~$g$ on an
oriented~$L^3$ and a real analytic section~$\text{I\!I}$ of
of the bundle~$S^3_0(T^*\!L)$, i.e., the traceless cubic forms
on~$L$, then there exists an isometric embedding of~$L$
as a special Lagrangian submanifold of a Calabi-Yau~$M^6$
with the given~$\text{I\!I}$ as second fundamental form.
The proof proceeds exactly as in Theorem~$1$, with the modification
that one uses the information in~$\text{I\!I}$ to choose
different restraining manifolds~$Z_3$, $Z_4$, and $Z_5$.
Details will be left to the interested reader.
\endremark

\remark{Special Lagrangian foliations}
Finally, if~$g$ is a real analytic metric on the $m$-torus~$T^m$ with
the property that every nonzero $g$@-harmonic 1-form on~$T^m$ 
is nowhere vanishing, then for any realization of~$(T^m,g)$
as a special Lagrangian manifold in a Calabi-Yau~$(M^{2m},\omega,\Upsilon)$,
the set of nearby special Lagrangian tori foliate~$M$ in
a neighborhood of~$T^m$.  (This is a consequence of McLean's
description of the moduli space.)  For a discussion of 
the significance of these foliations for string theory and mirror 
symmetry, see~\cite{SYZ}.  For a discussion of the mathematical aspects
of such foliations, see~\cite{Hi1,Hi2}.

It is easy to see that there are many metrics on the~$m$@-torus
for which a basis of the harmonic~$1$-forms is everywhere linearly
independent.  In fact, when~$m=2$, all the metrics on the
torus have this property.

On the other hand, suppose that~$m>2$ and
that one has such a metric~$g$ on~$T^m$.
One can identify~$T^m$ with~$\bbR^m/\bbZ^m$ in such a way 
that the differentials~$dx^i$ are a basis for 
the $g$@-harmonic 1-forms on~$T^m$ so suppose that
this has been done.%
\footnote{The essential ambiguity in this 
identification is an element of~$\SL(m,\bbZ)$.}
The metric~$g$ has the form~$g = g_{ij}\,dx^i\,dx^j$ 
for some functions~$g_{ij}=g_{ij}$ on the torus~$T$. 
Let~$\Delta = \det(g_{ij})$ and 
set~$h = \Delta^{1/2} (g_{ij})^{-1}$.  Then~$h = (h^{ij})$
is a symmetric, positive definite matrix on~$T^m$
and the condition that the forms~$dx^i$ be $g$@-harmonic
becomes the~$m$ linear, first order equations
$$
\frac{\partial h^{ij}}{\partial x^j} = 0.
$$
Conversely, if~$h$ is any positive definite solution on~$T^m$ to
these linear equations, setting
$$
(g_{ij}) = \det(h)^{1/(m-2)} h^{-1}
$$
defines a metric~$g = g_{ij}\,dx^i\,dx^j$
on~$T^m$ for which the~$dx^i$ are $g$@-harmonic.

Thus, the space of such metrics (even in the real analytic 
category) is a convex open set in a linear space.
\endremark

\head 2. Coassociative Realization  \endhead

\proclaim{Theorem 2}
Let~$(L^4,g)$ be a closed, real analytic, oriented Riemannian $4$-manifold
whose bundle of self-dual $2$-forms is trivial.  
Then there exists a $\G_2$@-manifold~$(N^7,\phi)$ and an
isometric embedding~$\iota:L\to N$ whose image is a coassociative
$4$-manifold in~$N$.  Furthermore, $(N^7,\phi)$ and $\iota$
can be chosen so that~$\iota(L)$ is the fixed point
locus of an anti~$\G_2$@-involution~$r:N\to N$.
\endproclaim

\demo{Proof}
First, I will examine the differential system~$\cI$ constructed in~\S0.4 for
the case~$G=\G_2$ and show that it is regularly presented.  Throughout this
proof, lower case latin indices will range from~1~to~7.

As in \S0.3,
let~$G=\G_2\subset\SO(7)$ be the subgroup that stabilizes the $3$@-form
$$
\phi_0 = dx^{567} - dx^5\w(dx^{12}+dx^{34}) - dx^6\w(dx^{13}+dx^{42})
           - dx^7\w(dx^{14}+dx^{23}).
$$
The ring~$\Lambda^*(\bbR^7)^G$ is generated by~$\phi_0$ and the $4$@-form
$$
\ast\phi_0 = dx^{1234}-dx^{67}\w (dx^{12}+dx^{34})-dx^{75}\w(dx^{13}+dx^{42})
           - dx^{56}\w(dx^{14}+dx^{23}).
$$
In particular,~$\G_2$ is admissible.

One can now compute the spaces~$\euh_k$ for~$0\le k\le 7$.
This computation is essentially the same as that in~\cite{Br2, Proposition~2},
so I will not go into detail.  One has~$\euh_k = M_7(\bbR)$ for~$k<3$,
while~$\euh_3$ is defined by the single equation~$x^5_3-x^6_2+x^7_1=0$.
Next,~$\euh_4\subset\euh_3$ is defined by four more equations
$$
x^1_1+x^2_2+x^3_3+x^4_4=
 x^5_2+x^6_3+x^7_4=x^5_1+x^6_4-x^7_3=x^5_4-x^6_1-x^7_2=0.
$$
Then~$\euh_5\subset\euh_4$ is defined by ten more equations:
$$
x^1_1+x^2_2-x^3_3-x^4_4=x^3_1+x^4_2+x^1_3+x^2_4=x^4_1-x^3_2-x^2_3+x^1_4=0
$$
and
$$
x^1_5-x^6_4+x^7_3 = x^2_5-x^6_3-x^7_4 = x^3_5+x^6_2-x^7_1 = 
x^4_5+x^6_1+x^7_2 = x^5_5 = 0
$$
$$
2x^6_5-x^4_1-x^3_2+x^2_3+x^1_4 = 2x^7_5+x^3_1-x^4_2-x^1_3+x^2_4 = 0.
$$
Finally, from general 
considerations~$\euh_6=\eug_2+M_{7,6}(\bbR)\simeq\bbR^{21}$ 
and~$\euh_7 = \eug_2\simeq\bbR^{14}$.

In particular, letting~$c_i$ denote the codimension 
of~$\euh_i$ in~$M_7(\bbR)$, it follows that
$$
(c_0,c_1,\ldots,c_7) = (0,0,0,1,5,15,28,35).
$$
Since
$$
c_0 +\cdots + c_6 = 49 = 7\cdot(21 - 14),
$$
it follows by Cartan's Test that~$\G_2$ is strongly admissible~
\cite{Br2, Proposition~1} and that~$\G_2$ is regularly presented.  In fact,
using the result~\cite{Br2, Proposition 2} that~$\G_2$ acts transitively 
on the space of oriented~$p$-planes in~$\bbR^7$ except for~$p=3$ and $p=4$,
it can be shown that any conjugate of~$\G_2$ in~$\SO(7)$
is regularly presented. 

Next, let~$M = L\times\bbR^3$ with linear coordinates~$y^1,y^2,y^3$ on the
second factor and let~$r:M\to M$ be the involution 
defined by~$r(p,y) = (p,-y)$.  For notational simplicity, I will identify
$L$ with~$L\times{\bold0}$, the fixed point set of~$r$.  Let~$\pi:F\to M$
be the coframe bundle and define an involution (also denoted~$r$) on~$F$
by the rule
$$
r(u) = - u\circ r'(x)
$$
for~$u\in F_x$.  Note that~$\pi\bigl(r(u)\bigr) = r\bigl(\pi(u)\bigr)$ 
and, tracing through the definitions, that~$r^*\eta = -\eta$.  
Since~$-\text{I}_7$ does not lie in~$\G_2$, the involution~$r:F\to F$
does not preserve~$\widehat{\phi_0}$.  In fact, 
$$
r^*\bigl(\widehat{\phi_0}\bigr) = {} - \widehat{\phi_0}
\qquad\text{while}\qquad
r^*\bigl(\widehat{\ast\phi_0}\bigr) =  \widehat{\ast\phi_0}\,.
$$
In particular, $r$ preserves~$\cI$ and hence its integral manifolds.
Since~$-\text{I}_7$ commutes with the elements of~$\G_2$, it
follows that~$r$ descends to a well-defined involution of~$S = F/\G_2$,
also to be denoted by~$r$.

Let~$\cI$ be the ideal on~$F$ or~$S$ generated 
by~$d\bigl(\widehat{\phi_0}\bigr)$ and~$d\bigl(\widehat{\ast\phi_0}\bigr)$.
Then, since~$\G_2$ is regularly presented, any~$E\in V_7(\cI,\pi)$ is the 
terminus of a regular flag, namely the 
canonical flag~$\{E_i\,\vrule\,0\le i\le 7\}$ as defined in~\S0.5.  
In fact, by the remark above that any conjugate of~$\G_2$ is regularly
presented, it follows that {\it all\/} of the elements~$E$ in 
either~$V_i(\cI,\pi)$ or~$V_i(\cI,\bar\pi)$ are regular for~$0\le i\le 7$,
with the codimension of~$H(E)$ being~$c_i$ as defined above.
 
Moreover, since~$\G_2$ has been shown to be strongly admissible, 
it follows that any section~$\sigma:U\to S$ over an open~$N\subset M$
is torsion-free and so induces the structure of a $\G_2$-manifold on~$N$.

Now let~$(L^4,g)$ be as in the hypotheses of Theorem~2.
By hypothesis,~$\Lambda^2_+(TL)$ is smoothly trivial.
Since~$g$ is closed, real analytic and Riemannian, Bochner's 
result that the real analytic forms are dense in the smooth
forms in the uniform topology can be applied to show that
$\Lambda^2_+(TL)$ can be real analytically trivialized.  
I.e., there exist three, real analytic
self-dual 2-forms~$\Omega_1$, $\Omega_2$, and $\Omega_3$ on~$L$ so that
$$
\Omega_i\w\Omega_j = 2\delta_{ij}\ {\ast_g}1.
$$
It is an elementary result in linear algebra that, at every point~$x\in L$,
there exists an oriented $g$@-orthonormal basis~$\alpha_0,\ldots,\alpha_3$
of~$T_x^*L$ so that
$$
\align
\Omega_1(x) &= \alpha_0\w\alpha_1{+}\alpha_2\w\alpha_3\,,\\
\Omega_2(x) &= \alpha_0\w\alpha_2{+}\alpha_3\w\alpha_1\,,\\
\Omega_3(x) &= \alpha_0\w\alpha_3{+}\alpha_1\w\alpha_2\,.\\
\endalign
$$
($L$ may not be parallelizable, so there may not be 
a global $g$@-orthonormal coframing with this property.)

Consider the $3$-form on~$M$
$$
\varphi = dy^1\w dy^2\w dy^3-dy^1\w\Omega_1-dy^2\w\Omega_2-dy^3\w\Omega_3\,.
$$
By the linear algebra result stated above and the definition of~$\phi_0$
given in~\S0.4, at every point~$x\in M$ there exists
a linear isomorphism~$u:T_xM\to\bbR^7$ so that~$u^*(\phi_0) = \varphi_x$.
Consequently, $\varphi$ defines a $\G_2$@-structure on~$M$ and hence
corresponds to a section~$\bar\sigma:M\to S$ that 
satisfies~$\sigma^*\widehat{\phi_0}=\varphi$.  
The metric underlying this~$\G_2$@-structure is
$\bar h = g+(dy^1)^2+(dy^2)^2+(dy^3)^2$, the Hodge dual of~$\varphi$
with respect to~$\bar h$ is
$$
\ast_{\bar h}\varphi
= \ast_g1-dy^2\w dy^3\w\Omega_1-dy^3\w dy^1\w\Omega_2-dy^1\w dy^2\w\Omega_3\,, 
$$
and evidently~$\sigma^*\widehat{\ast\phi_0} = \ast_{\bar h}\varphi$.  

Let~$X_4 = \sigma(L)$.  Since~$d\varphi$
vanishes when pulled back to~$L$, it follows that~$X_4$
is an integral manifold of~$\cI$ that is transverse to~$\bar\pi$.
Moreover, since~$r^*\varphi = -\varphi$, it follows that~$X_4\subset S$ 
is pointwise fixed under~$r$.  Since~$X_4$ is ~$\bar\pi$@-transverse,
it is a regular integral manifold of~$\cI$.

 From this point, the proof of Theorem~2 will follow the same lines as
the proof of Theorem~1.  There are some details to check, since the
construction of the `restraining' manifolds needed for the application
of the Cartan-K\"ahler theorem requires some care in the absence of
a global parallelization of~$L^4$, but this is a detail best 
left to the reader.  The crucial point is that one can define
subspaces~$W_5\subset W_{15}\subset W_{28}\subset M_7(\bbR)$ 
of dimensions~$5$, $15$, and~$28$ respectively that satisfy
$$
W_5\cap\euh_4 = W_{15}\cap\euh_5 = W_{28}\cap\euh_6 = (0),
$$
are invariant under conjugation by the element
$$
R = \pmatrix \text{I}_4 & 0 \\ 0 & -\text{I}_3\\ \endpmatrix,
$$
and are invariant under the subgroup~$\SU(2)\subset\G_2$ consisting
of the transformations that fix~$\phi_0$,~$dx^5$, $dx^6$, and $dx^7$.  
These spaces are then used to construct the restraining manifolds 
$Z_4$, $Z_5$, and $Z_6$ below, in a manner completely analogous to 
the constructions in the proof of Theorem~1.  (The $\SU(2)$@-invariance 
of the~$W_k$ allows one to define these restraining manifolds without 
reference to a coframing on~$L$.)

The rest of the proof can now be described as follows:

First, one constructs an $r$@-invariant real analytic manifold
$Z_4\subset S$ that contains~$X_4$, submerses
onto~$L\times\bbR^1\subset M$ with fibers of dimension~$5$,
and satisfies~$\dim\bigl(T_xZ_4\cap H(T_xX_4)\bigr) = 5$
for all~$x\in X_4$.
(That this is the appropriate dimension for the fibers follows from the 
calculation of~$\euh_4$.) Applying the Cartan-K\"ahler Theorem produces an
integral manifold~$Y_5$ of~$\cI$ that satisfies~$X_4\subset Y_5\subset Z_4$.
The intersection~$X_5 = Y_5\cap r(Y_5)$ is $r$-invariant and, by the
uniqueness part of the Cartan-K\"ahler theorem, it is a $5$@-dimensional
integral manifold of~$\cI$.
By shrinking~$X_5$ if necessary, one can ensure that it is the 
graph of a section of~$S$ over an open neighborhood~$N_5$ of~$L$ 
in~$L\times\bbR^1\subset M$.  In particular,~$X_5$ is 
$\bar\pi$@-transverse and so must be a regular integral manifold of~$\cI$.

Next, one constructs an $r$@-invariant real analytic manifold
$Z_5\subset S$ that contains~$X_5$, submerses onto an open neighborhood 
of~$L$ in~$L\times\bbR^2$ with fibers of dimension~$15$,
and satisfies~$\dim\bigl(T_xZ_5\cap H(T_xX_5)\bigr) = 6$ for all~$x\in X_5$.
(That this is the appropriate dimension for the fibers follows from the
calculation of~$\euh_5$.) Applying the Cartan-K\"ahler Theorem constructs an
integral manifold~$Y_6$ of~$\cI$ that satisfies~$X_5\subset Y_6\subset Z_5$.
The intersection~$X_6 = Y_6\cap r(Y_6)$ is $r$-invariant and, by the
uniqueness part of the Cartan-K\"ahler theorem, it is a $6$@-dimensional
integral manifold of~$\cI$.
By shrinking~$X_6$ if necessary, one can ensure that it is
is the graph of a section of~$S$ over an open neighborhood~$N_6$ of~$L$ 
in~$L\times\bbR^2\subset M$. In particular,~$X_6$ is $\bar\pi$@-transverse 
and so must be a regular integral manifold of~$\cI$.

Finally, one constructs an $r$@-invariant real analytic manifold
$Z_6\subset S$ that contains~$X_6$, submerses onto an open neighborhood 
of~$L$ in~$L\times\bbR^3$ with fibers of dimension~$28$,
and satisfies~$\dim\bigl(T_xZ_6\cap H(T_xX_6)\bigr) = 6$ for all~$x\in X_6$.
(That this is the appropriate dimension for the fibers follows from the
calculation of~$\euh_6$.) Applying the Cartan-K\"ahler Theorem constructs an
integral manifold~$Y_7$ of~$\cI$ that satisfies~$X_6\subset Y_7\subset Z_6$.
The intersection~$X_7 = Y_7\cap r(Y_7)$ is $r$-invariant and, by the
uniqueness part of the Cartan-K\"ahler theorem, it is a $7$@-dimensional
integral manifold of~$\cI$.
By shrinking~$X_7$ if necessary, one can ensure that it is the graph of a 
section~$\sigma$ of~$S$ over an open neighborhood~$N=N_7$ of~$L$ 
in~$M = L\times\bbR^3$.

Let~$\phi=\sigma^*(\widehat{\phi_0})$. 
By construction,~$d\phi=d\,{\ast}\phi=0$ and $r^*\phi = -\phi$.  
Thus, $(N,\phi)$ is a $\G_2$@-manifold with anti@-$\G_2$ involution~$r$.  
Now, $\phi$ and $\varphi$ agree along~$L$, the fixed locus of~$r$.  
Thus,~$L$ is an integral manifold 
of~$\phi$ since it is visibly an integral manifold of~$\varphi$.  
Moreover, if~$h$ is the metric on~$N$ associated to~$\phi$, 
then $h$ and $\bar h$ agree along~$L$.  

Consequently, the inclusion of~$L$ into~$N$ is an isometric embedding 
of~$(L,g)$ into~$(N,h)$ as a coassociative submanifold, in particular,
as the fixed locus of the anti@-$\G_2$ involution~$r$.\qed
\enddemo

\remark{Topological conditions}
If~$L$ is an oriented Riemannian $4$-manifold and~$\Lambda^2_+(TL)$ 
is trivial, then the structure group of~$TL$ 
can be reduced to a copy of~$\SU(2)\subset\SO(4)$.  In particular,
$L$ must be spin.  Also, the first Pontrjagin class of~$\Lambda^2_+(TL)$
must vanish.  This class is~$p_1(TL)+2e(TL)$, so the signature
theorem~\cite{MS} implies
$$
b_2^-(L) + 4 b_1(L) = 5b_2^+(L) + 4.
$$
Conversely, if an oriented, spin $4$-manifold~$L$ satisfies this relation,
then for any metric~$g$, its bundle of self-dual 2-forms is 
topologically trivial.  Examples of such compact manifolds 
are the 4-torus and the K3 surface, but these are not the only possibilities, 
of course.

The reader may wonder whether or not the triviality of $\Lambda^2_+(TL)$ 
is really necessary for the conclusion of Theorem~2. Certainly the proof
via Cartan-K\"ahler theory requires it, but there could conceivably be
a way around this.  However, a careful look at the construction in the
proof of Theorem~2 shows that the choice of trivialization 
of~$\Lambda^2_+(TL)$ does affect the resulting $\G_2$@-structure on~$N$.  
I do not know how one could `patch' to avoid this.  Nevertheless,
there are known compact examples~\cite{BSa, Jo} for which~$\Lambda^2_+(TL)$
is not trivial.

Finally, the reader may have noticed that the only place that the
compactness of~$L$ was used was to prove that the topological triviality
of~$\Lambda^2_+(TL)$ implies that this bundle is real analytically 
trivial.  Naturally, one would expect this to be true even without
the compactness hypothesis.  I have been unable to find an explicit 
proof of this in the literature, but perhaps the reference~\cite{Sh}
is relevant.
\endremark

\remark{Coassociative fibrations}
The analog of special Lagrangian fibrations in the coassociative
setting is, of course, coassociative fibrations.  If a~$\G_2$@-manifold~
$(N,\phi)$ has a fibration~$\beta:N\to B^3$ whose fibers are 
compact coassociative submanifolds, then the normal bundle of each
fiber~$L = L_b = \beta^{-1}(b)$ is trivial.  Since this normal bundle is 
canonically isomorphic to~$\Lambda^2_+(TL)$, this latter bundle is
trivial.  Moreover, by McLean's description of the moduli space
of coassociative submanifolds, it follows that there exist three
harmonic, self-dual 2-forms~$\Omega_1$, $\Omega_2$, $\Omega_3$
in~$\Omega^2_+(L)$ that are everywhere linearly independent.  

Of course, these $2$@-forms are closed and satisfy
$$
\Omega_i\w\Omega_j = 2a_{ij}\,\,{\ast1}
$$
for some~$a = (a_{ij})$ that is symmetric and positive definite at
each point of~$L$.

Conversely, suppose that~$L^4$ is a compact, real analytic 
$4$@-manifold and that there are three closed, real analytic $2$@-forms 
$\Omega_1$, $\Omega_2$, and $\Omega_3$ on~$L$ that satisfy the 
nondegeneracy condition that 
$$
\Omega_i\w\Omega_j = 2a_{ij}\,\,\Phi
$$
for some~$a = (a_{ij})$ that is symmetric and positive definite at
each point of~$L$ and some nonvanishing $4$@-form~$\Phi$ on~$L$.

Then there is a unique orientation and real analytic conformal 
structure on~$L$ so that the forms~$\Omega_i$ are a basis for the
self-dual $2$-forms on~$L$.  Choosing a real analytic metric~$g$
in this conformal class, one can apply Theorem~$2$ and McLean's
description of coassociative moduli to construct a%
\footnote{in fact, infinitely many} 
$\G_2$@-manifold~$(N^7,\phi)$ so that~$(L,g)$ is isometrically embedded 
as the $\bold0$@-fiber of a coassociative fibration~$\beta:N^7\to B^3$, 
where~$B^3$ is a neighborhood of~$\bold0$ in~$\bbR^3$.

Now, the obvious examples of such triples~$(\Omega_1,\Omega_2,\Omega_3)$
are the parallel examples on the flat 4-torus~$T^4$ 
and the parallel $2$@-forms on a K3 surface endowed with its Calabi-Yau
metric.  However, these are by no means the only ones.  
In fact, given one such triple~$(\bar\Omega_1,\bar\Omega_2,\bar\Omega_3)$,
one can construct many others in the form
$$
(\Omega_1,\Omega_2,\Omega_3) 
=(\bar\Omega_1{+}d\alpha_1,\bar\Omega_2{+}d\alpha_2,\bar\Omega_3{+}d\alpha_3)
$$
where the $\alpha_i$ are real analytic $1$-forms whose exterior derivatives
are small in the uniform norm.  (Generically, these will not be parallel
with respect to any metric on~$L$.)  
\endremark

\remark{Prescribed second fundamental forms}  The construction in
Theorem~$2$ produces a $\G_2$-manifold in which~$L$ appears as
a totally geodesic submanifold.  However, modifying the construction
allows one to get other second fundamental forms.  

For a coassociative submanifold~$L\subset N$ where~$(N,\phi)$ is
a $\G_2$@-manifold, the second fundamental form~$\text{I\!I}$ takes
values in a certain 15-dimensional subbundle~$S\subset\Lambda^2_+(TL)\otimes
S^2_0(T^*\!L)$  that can be defined using only the metric~$g$ on~$L$ and
not the ambient metric or connection on~$N$. (The bundle~$S$ is
the coassociative analog of~$S^3_0(T^*\!L)$, which appeared as the 
receiving bundle for second fundamental forms of special Lagrangian
submanifolds.)  When~$L$ is parallelizable and a real analytic 
section~$\text{I\!I}$ is specified, one can show that 
the restraining manifolds~$Z_i$ used in the 
Cartan-K\"ahler construction can be chosen so as to produce a 
$\G_2$-manifold into which~$(L,g)$ embeds isometrically as
a coassociative submanifold with second fundamental form equal to~
$\text{I\!I}$.  Details are left to the reader.
\endremark

\enddocument